\documentclass[12pt]{article}
\usepackage[ruled,vlined]{algorithm2e}
\usepackage{amsmath}
\usepackage{amsthm}
\usepackage{amssymb}
\usepackage{color}
\usepackage{subcaption}
\usepackage[page]{appendix}
\usepackage{mathtools}
\usepackage{caption}
\usepackage{bm}
\usepackage{hyperref}

\usepackage{booktabs}
\usepackage{graphicx}
\usepackage{authblk}
\usepackage{float}

\usepackage[margin=1in]{geometry}

\usepackage{etoolbox}
\AtEndEnvironment{section}{\setcounter{claim}{0}}
\theoremstyle{definition}

\theoremstyle{remark}

\DeclarePairedDelimiterX{\infdivx}[2]{(}{)}{%
  #1\;\delimsize\|\;#2%
}

\title{The Inverse Micromechanics Problem given Dielectric Constants for Isotropic Composites with Spherical Inclusions}

\date{\today}

\begin{document}
\author[1]{Athindra Pavan}
\author[2]{Swaroop Darbha}
\author[1,3]{Bj\"{o}rn Birgisson}
\affil[1]{School of Environmental, Civil, Agricultural and Mechanical Engineering, College of Engineering, University of Georgia, 200 DW Brooks Drive, Athens, GA 30602, USA}
\affil[2]{Department of Mechanical Engineering, Texas A\&M University, College Station, TX, USA}
\affil[3]{Office of the Provost, Nazarbayev University, 53, Kabanbay batyr Ave., 010000, Astana, Kazakhstan}

\maketitle
\begin{abstract}
    In this article, convex optimization is introduced as a promising tool to study Eshelby based inverse micromechanics problems. The focus is on inverse micromechanics using the Eshelby-Mori-Tanaka model given the dielectric constants of the composite material and of all of its components. The model is exactly the same for the conductivity properties as well. This choice of model is made since the model is fairly simple and has a closed form analytical solution for the case of spheroidal inclusions as well. The forward or direct micromechanics problem deals with the determination of effective properties of a composite material given the properties of its components and microstructural information. The focus is on isotropic composites and the distribution of inclusions is assumed to be such that this holds. The inverse micromechanics problem considered in this paper deals with the determination of microstructural information given the properties of the composite material and all of its components. Since in this paper, isotropy of the composite and only spherical inclusions are considered, the goal is to determine just volume fractions of the components of the composite material. The inverse problem is formulated as a Linear Programming problem and is solved. Before this, the inverse problem and certain important variants of it are examined through the lens of convex optimization. Lastly, promising results regarding the relationship between dispersive materials, noise in measurements, and quality of obtained volumetric splits are showcased. The scope of the use of convex optimization in inverse micromechanics is discussed. 
\end{abstract}

\setcounter{section}{-1} 
\section{Notation}
\begin{itemize}
    \item Lowercase alphbets (Greek or otherwise) in boldface are vectors. The main vector in this paper, $\bm{\phi}$, for example.
    \item Scalars are labeled using subscript as $t_k$, etc. Vectors and matrices are labeled using superscript so that the subscript may be available to refer to entries of the object. So $\bm{\phi^S}$ is a vector, and $\phi^S_i$ is the $i^{\text{th}}$ entry of the vector.
    \item $\bm{1^n}$ is a vector comprising $n$ 1s and $\bm{I^{n}}$ is an identity matrix in $\mathbb{R}^{n\times n}$. The notation used in Boyd's textbook \cite{boyd2004convex} of $\left(1^\top a+b\right)$ to denote the sum of a dot product of a vector of 1s with another vector `a' and scalar `b' is not used. The use of non-boldface characters for vectors and scalars is ambiguous. So, instead, one has $\left(\bm{1^n}\cdot\bm{a}+b\right)$. 
\end{itemize}

\section{Introduction}
Micromechanics has been a subject of study since as
early as 1880 undertaken by many researchers like Ernst
Gustav Kirsch. The starting point for this paper is Maxwell's work \cite{maxwell1873treatise} on the uniform internal field of an ellipsoidal body in a homogeneous medium under a uniform external field. Stratton's textbook \cite{stratton2007electromagnetic} first published in 1941 presents the same results in a manner that is more directly applicable to the study of composite materials. These results are used for the dielectric permittivity micromechanics problem. The analog of this result for the linear elasticity micromechanics problem is Eshelby’s result \cite{Esh1957} showing that the stress and strain fields within an ellispoidal inhomogeneity embedded in an infinite matrix are uniform when the system is loaded at infinity. These results have since been used to develop many homogenization (or mean-field averaging or micromechanics) models. Examples of such models are the Mori-Tanaka (M-T) model, and the Self-Consistent Scheme (SCS) among others (as described in great detail in \cite{mura2013micromechanics}, \cite{qu2006fundamentals}, \cite{torquato2002random}, etc.). Since in the micromechanics community this result is often called the Eshelby-Mori-Tanaka or the Mori-Tanaka-Benveniste model, this name will be used. The focus of this study is inverse micromechanics using the Eshelby-Mori-Tanaka model. For the dielectric permittivity problem, the formulation of the Eshelby-Mori-Tanaka model presented in \cite{benveniste1986effective} for thermal conductivity is used. The same problem formulation is valid for the dielectric permittivity problem as well (see \cite{torquato2002random}). The direct or forward problem deals with determining the effective properties of the composite material given properties of its components and microstructural information (volume fractions, inclusion shapes and orientations, probability density function, etc). The inverse problem studied in this work deals with determining the microstructural information given the effective properties of the composite and the properties of all components. Other variations of the inverse problem have been studied.

The composite considered in this work is such that all its components are linear materials (properties are linear) that are isotropic and homogeneous. In addition to this, the nature of the microstructure of the composite is such that it comprises many spherical inclusions distributed in such a manner that the composite as a whole can be studied as an isotropic continuum with linear behavior. This is called statistical isotropy and is discussed in detail in \cite{hashin1983analysis}. If this assumption is not made, the Probability Density Functions (PDF) of inclusion materials must be specified in the study. This is an added layer of complexity that is not addressed in this paper. The Eshelby-Mori-Tanaka model, which is used in this paper, for the general composite material is known to have realizability issues as discussed in \cite{Norris1989}, \cite{BERRYMAN1996149}, etc. Realizability is defined in 2 ways: (1) the effective properties abide by the Hashin-Shtrikman bounds \cite{HASHIN1963127}, or (2) the volumetrics considered correspond to a well defined real microsctructure. For the case considered here (isotropic composite material with spherical inclusions), these are not issues. Specifically because the focus is on inverse micromechanics with real dielectric constant measurements driving the determination of volumetrics.

Inverse micromechanics problems and the concept of an equivalent microstructure are related and quite challenging to study. Several versions of the inverse micromechanics problem can be studied. In \cite{ogierman2018inverse} for example, the microstructure and volume fraction of inhomogeneities are known and the properties of the component materials are determined using an evolutionary algorithm. \cite{LYDZBA201820} studies the inverse thermal conductivity micromechanics problem and solves for a Probability Density Function (PDF) that describes the microstructure of an isotropic random heterogeneous material whose effective properties and component material properties are known. Infinite families of oblate spheroidal shapes with semi-axis aspect ratio considered to be in the interval (0,1) are assumed in \cite{LYDZBA201820}. In \cite{lydzba2019principle} they have studied equivalent microstructures comprising only oblate spheroidal inclusions for isotropic composites. They show that the equivalent microstructure is independent of the properties of inhomogeneities but depends on the volume fractions of components. In \cite{sevostianov2012concept}, it is seen that an ``average ellipsoidal shape" can be used to replace diverse ellipsoidal inclusions (varying aspect ratios) in an isotropic composite material. This means that a microstructure with randomly oriented identical ellipsoidal inclusions of this ``average shape" will have the same effective properties as the original microstructure with diverse ellipsoidal inclusions. In \cite{sevostianov2012concept}, it is also shown that it is possible to have an equivalent microstructure for an isotropic two-component composite material with inclusions of arbitrary shape to be replaced by an equivalent microstructure having randomly oriented identical ellipsoidal inclusions. In determining microstructural parameters to solve a variety of inverse problems, simulated annealing (SA) has been widely used (\cite{vcapek2009stochastic}, \cite{PhysRevE.59.5596}, \cite{rintoul1997reconstruction}, \cite{rozanski2011digital}, \cite{talukdar2002reconstruction}, \cite{talukdar2002stochastic}, \cite{yeong1998reconstructing}, etc.). These evolutionary algorithms use a lot of computational power. Convex optimization is a powerful mathematical tool that is used to solve optimization problems much faster than purely computational methods like simulated annealing, differential evolution, etc. However, a certain structure is required of the problem in order to use this tool. This work shows certain problem formulations that have the desired structure to use this powerful tool. Analytical papers like \cite{kachanov2005quantitative}, \cite{sevostianov2012concept}, \cite{lydzba2019principle}, that have been cited earlier in this paragraph are few in number and there is scope for much more work in the space of inverse micromechanics. 

There has been extensive work in the field of inverse micromechanics but not much using convex optimization to formulate optimization problems to solve for volumetric splits as is done in this paper. \cite{milton2012bounds}, \cite{kang2013bounds}, and \cite{thaler2013bounds} are examples of the use classical micromehcanics methods such as analysis of PDEs and the translation method (see \cite{cherkaev2012variational}) to arrive at bounds on volume fractions of inclusions based on known information of material properties (of the composite and its components). These methods do not use the shape information given by the Eshelby result \cite{Esh1957} which is leveraged in this paper.

Convex optimization is the parent branch of mathematics to which commonly studied problems like the linear programming problem, the quadratic programming problem, the least squares problem, etc, belong. The goal of the subject is to characterize and solve optimization problems. The hardest challenge one faces when using convex optimization is in trying to show that the problem at hand is indeed a convex optimization problem. However, for the case of the Eshelby-Mori-Tanaka model with spherical inclusions and isotropy, it is very simple. The reader may refer to \cite{boyd2004convex} for an in-depth explanation of concepts in convex optimization. 

To the best of the authors' knowledge, convex optimization has not been used in the field of Eshelby based inverse micromechanics in the manner in which it is introduced in this paper. Convex optimization is used here to formulate an inverse micromechanics problem utilizing the Eshelby-Mori-Tanaka mean-field averaging model for the complex dielectric constant of a composite material system. The use of the Eshelby result leverages shape based information. This work differs from classical micromechanics in that the actual volumetric split is the goal. Not bounds as provided in \cite{thaler2013bounds}, \cite{kang2013bounds}, and \cite{milton2012bounds}. The final problem formulation used in this work given complex dielectric constant measurements of the composite material and its components can potentially be used in the determination of the volumetric split at every location of a composite material. The problem strictly works only for composites with spherical inclusions. However, with the use of equivalent microstructures as discussed in Section \ref{section: Equivalent Microstructures}, extensions should be possible. Considering more sophisticated microstructures with the general ellipsoid, more generally applicable and better performing optimization problems can also be formulated. Currently, the admissible noise level is low, but this work highlights the role of the dispersion strength of the components of a material system in determining the volumetric split. The problem can be modified to include other physical information that is known about the system being studied so as to reduce the dependence on the dispersion strength of the components.  

In section \ref{section: theoretical Bkgd}, a brief theoretical background for the forward problem is given. In section \ref{section: Optimization problem}, the inverse problem is formulated as an optimization problem. In section \ref{section: characterization}, the optimization problem is characterized. Section \ref{section: solutions} deals with analytical solutions to a couple of variations of the formulated optimization problem. Section \ref{section: Final 2 formulations} discusses two formulations of the optimization problem of which \ref{subsection: Final Formulation} is the final optimization problem formulated for the class of composite material systems studied in this paper. Section \ref{section: application} deals with applications of the proposed optimization problem to three different material systems. Since only dielectric constants of all component materials and of the composite material itself are used, the applications show the need for a strong dispersive component in the composite material system to ensure the predicted volumetric split to be accurate. If dielectric constant information is not available, other physical constraints that are available can be included. The requirements of the rank and minimum singular value of the constraint sensitivity matrix ($\bm{G}$) is discussed in this section to aid researchers and engineers formulate a suitable optimization problem for the application at hand. Section \ref{section: Equivalent Microstructures} contains a short note on equivalent microstructures for material systems with complex microstructures. The paper concludes with Section \ref{section: conclusion} which discusses all the results obtained.

\section{Theoretical Background}\label{section: theoretical Bkgd}

The Eshelby-Mori-Tanaka model to compute the effective dielectric constant ($\tilde{\bm{\varepsilon}}$) for an $n$-component composite material with isotropic components and randomly oriented spheroidal inclusions is as follows:

\begin{align}\label{eqn: MTB}
\begin{split}
    \tilde{\bm{\varepsilon}}=\frac{\sum_{i=0}^{n-1}\left(\phi_i\varepsilon_i\tilde{\bm{R}^i}\right)}{\sum_{i=0}^{n-1}\left(\phi_i\tilde{\bm{R}^i}\right)};\quad R_i =\frac{1}{3}\sum\limits_{j=1}^{3}\left[1+A_{jj}^i\frac{\varepsilon_i-\varepsilon_0}{\varepsilon_0}\right]^{-1};\quad \tilde{\bm{R}^i}=R_i\bm{I}
\end{split}
\end{align}
Where in the above, associated with the $i^{\text{th}}$ inclusion material, are the following:
\begin{itemize}
    \item $\bm{R}^i$ is a concentration tensor.
    \item $\bm{A}^i$ is the symmetric depolarization tensor (analog to the Eshleby tensor).
    \item $\phi_i$ is the volume fraction.
\end{itemize}  
$\bm{I}$ is the second order identity tensor. The tilde $\tilde{\cdot}$ is used to represent an average over all orientations. $\bm{A}^i$ is given by: 
\begin{align}\label{eqn: Main equations}
\begin{split}
    &\bm{A}^i=\begin{bmatrix}
        Q_i & 0 & 0\\
        0 & Q_i & 0\\
        0 & 0 & 1-2Q_i
    \end{bmatrix};\quad \text{where,}\\
    Q_i&=\frac{\alpha_i}{2(\alpha_i^2-1)^{3/2}}(\alpha_i(\alpha_i^2-1)^{1/2}-\cosh^{-1}(\alpha_i))\quad \text{Prolate shape, $\alpha_i>1$}\\
    &=\frac{\alpha_i}{2(1-\alpha_i^2)^{3/2}}(\cos^{-1}(\alpha_i)-\alpha_i(1-\alpha_i^2)^{1/2})\quad \text{Oblate shape, $\alpha_i<1$}\\
    &\text{Where $\alpha_i=l_i/d_i$ is the aspect ratio of the spheroidal inclusion}
\end{split}
\end{align}

In this paper, only spherical inclusions are considered. So, $\bm{A}^i=\frac{1}{3}\bm{I}$ and $\tilde{\bm{R}}^i=\frac{3}{3+r_i}\bm{I}$ where $r_i=\frac{\varepsilon_i}{\varepsilon_0}-1$. For more details on the forward model, refer to \cite{benveniste1986effective}.

\section{Optimization Problem}\label{section: Optimization problem}
The inverse problem involves solving for the microstructural parameters of the composite material given the effective properties of the composite and the properties of its components. In this work, a simplified problem is solved. The material properties of the composite material are assumed to be isotropic due statistical isotropy (see \cite{hashin1983analysis}) as a result of the manner in which the inclusions are distributed within the composite material. Given a composite material, a couple of preliminary experiments will help determine whether the material is isotropic or not. After this, the concept of an equivalent microstructure could also be used (see Section \ref{section: Equivalent Microstructures}). If the microstructure is known to be such that the composite material system is isotropic and all inclusions are spherical in shape, then the optimization problem developed here can be directly applied. For the general ellipsoidal inclusion, more work needs to be done, and this is an effort toward that endeavor.   
The optimization problem in question for the current inverse Eshelby-Mori-Tanaka mean-field averaging problem is as follows:
\begin{equation}\label{eqn: main optimization problem}
\boxed{
\begin{aligned}
    \underset{\bm{\phi}}{\text{Minimize}}\> &\eta\coloneqq \left\lvert\overline{\varepsilon}-\hat{\varepsilon}\right\rvert\\
    \text{Subject to: }&\\
    &\phi_i\in(0,1)\>\forall\> i\in\{0,1,2,\cdots,(n-1)\};\\
    &\bm{1^n}\cdot\bm{\phi}=1
\end{aligned}}
\end{equation}
Here, $\overline{\varepsilon}=\tilde\varepsilon/\varepsilon_0$ and $\tilde\varepsilon$ is computed using the method described in the previous section. $\hat{\varepsilon}=\hat{\hat\varepsilon}/\varepsilon_0$ and $\hat{\hat\varepsilon}$ is the experimentally determined dielectric constant of the composite material. In this article, the following results are shown for the optimization problem \ref{eqn: main optimization problem}:
\begin{itemize}
    \item It is a quasiconvex optimization problem.
    \item It is strictly quasiconvex for $n=2$. 
    \item The set of all minimizers $\mathcal{M}$ is determined analytically for $n\ge3$.    
\end{itemize} 

\section{Characterization of the Optimization Problem}\label{section: characterization}

\subsection{\texorpdfstring{The Eshelby-Mori-Tanaka Model is Quasilinear in $\bm\phi$}{The Eshelby-Mori-Tanaka Model is Quasilinear in phi}}
Equation \ref{eqn: MTB} reduces to the following for spherical inclusions after non-dimensionalizing by dividing by $\varepsilon_0$:
\begin{eqnarray}\label{eqn: main spherical in terms of r}
    \overline\varepsilon=\frac{\phi_0+\sum\limits_{i=1}^{n-1}\left(\frac{3(1+r_i)\phi_i}{(3+r_i)}\right)}{\phi_0+\sum\limits_{i=1}^{n-1}\left(\frac{3\phi_i}{(3+r_i)}\right)}; \text{ where, }r_i=\frac{\varepsilon_i}{\varepsilon_0}-1; r_i\in(-1,0)\cup(0,\infty).
\end{eqnarray}
Simplifying further to make analysis easier,
\begin{eqnarray}
    \overline\varepsilon&=&\frac{\bm{p}\cdot\bm{\phi}}{\bm{q}\cdot\bm{\phi}};\quad 1_{\{i=0\}}\coloneqq\begin{cases}
        0\> \forall\> i>0,\\
        1,\>i=0
    \end{cases};\quad i\in\{0,1,2,\cdots,(n-1)\}\nonumber\\
    p_i&\coloneqq& \left(1_{\{i=0\}}\right)\prod\limits_{k=1}^{n-1}(3+r_k)+3\left(1-1_{\{i=0\}}\right)\left((1+r_i)\prod\limits_{j\neq i,j=1}^{n-1}(3+r_j)\right)\label{eqn: N}\\
    q_i&\coloneqq&\left(1_{\{i=0\}}\right)\prod\limits_{i=1}^{n-1}(3+r_i)+3\left(1-1_{\{i=0\}}\right)\left(\prod\limits_{j\neq i,j=1}^{n-1}(3+r_j)\right)\label{eqn: D}
\end{eqnarray}
The function $\overline\varepsilon$ is a linear fractional function (see Boyd \cite{boyd2004convex}). This type of function is known to be quasilinear when $\bm{q}\cdot\bm\phi>0$ for all $r_i$ and $\bm\phi$. Thus $\bm{q}\cdot\bm{\phi}$, recalling that $r_i\in(-1,0)\cup(0,\infty)$, is:
\begin{gather}
\bm{q}\cdot\bm\phi=\left(1_{\{i=0\}}\right)\left(\prod_{i=1}^{n-1}\underbrace{(3+r_i)}_{>0}\right)\underbrace{\phi_0}_{>0}+3\left(1-1_{\{i=0\}}\right)\underbrace{\sum_{i=1}^{n-1}\left(\prod_{j\neq i;j=1}^{n-1}(3+r_j)\phi_i\right)}_{>0}>0
\end{gather}

\subsection{Strict Quasiconvexity of the Inverse Problem for 2 Component Composite Material with Spherical Inclusions}

For the case of 2 component composite materials, problem \ref{eqn: main optimization problem} is shown to be strictly quasiconvex. In this case, since a global minimum exists, this means that a unique minimizer $\left(\phi_0^\star ,\phi_1^\star \right)$ exists. This matches a result we have from micromechanics literature (see \cite{kachanov2005quantitative} Figure 3). To show this, first note:
\begin{equation}\label{eqn: ehat in range n=2}
    \min\left(1,\frac{\varepsilon_1}{\varepsilon_0}\right)<\hat\varepsilon<\max\left(1,\frac{\varepsilon_1}{\varepsilon_0}\right);\quad \hat{\varepsilon}\in\mathrm{Im}\left(\overline{\varepsilon}\right)
\end{equation}
Since $\hat\varepsilon$ is in the image of $\overline{\varepsilon}$, and since $\overline{\varepsilon}$ is continuous, it suffices to show that $\overline{\varepsilon}$ is monotone in $\phi_1$. Restating $\overline{\varepsilon}$ (equation \ref{eqn: main spherical in terms of r}) in terms of only $\phi_1$ for $n=2$ and finding its derivative:
\begin{equation}
    \overline{\varepsilon}=\frac{3+r_1+2r_1\phi_1}{3+r_1-r_1\phi_1};\quad \overline{\varepsilon}_{,\phi_1}=\frac{3r_1\left(3+r_1\right)}{\left(3+r_1-r_1\phi_1\right)^2} 
\end{equation}
Thus, for a given $r_1$, $\overline{\varepsilon}$ is monotone in $\phi_1$ and thus problem \ref{eqn: main optimization problem} is a strictly quasiconvex optimization problem. 

\subsection{\texorpdfstring{Quasiconvexity for $n$ Component Composite Material with Spherical Inclusions}{Strict Quasiconvexity for n Component Composite Material with Spherical Inclusions}}\label{subsec: SQC phi}

Before proving the main results of the current work within this subsection, some more definitions will be made here. Restating the optimization problem,

\begin{equation}\label{eqn: main optimization problem full definition}
\boxed{
\begin{alignedat}{3}
\underset{\bm{\phi},\eta}{\text{Minimize: }}&\eta\coloneqq\left\vert\frac{\bm{p}\cdot\bm{\phi}}{\bm{q}\cdot\bm{\phi}}-\hat{\varepsilon}\right\vert=\left\lvert\frac{\bm{u}\cdot\bm{\phi}}{\bm{q}\cdot\bm{\phi}}\right\rvert\\
\text{Subject to: }\\
&\bm{\phi}\succeq\bm{0};\quad\bm{1^n}\cdot\bm{\phi}=1& (C1)-(C2)\\
\hline
    1_{\{i=0\}}&\coloneqq\begin{cases}
        0,i>0,\\
        1,i=0
    \end{cases};\quad r_i\coloneqq\frac{\varepsilon_i}{\varepsilon_0}-1\>\forall i\in\{1,\cdots,(n-1)\};\\
    \hline
    p_i\coloneqq \left(1_{\{i=0\}}\right)&\prod\limits_{j=1}^{n-1}(3+r_j)+3\left(1-1_{\{i=0\}}\right)\left((1+r_i)\prod\limits_{j\neq i,j=1}^{n-1}(3+r_j)\right);\\
    q_i\coloneqq\left(1_{\{i=0\}}\right)&\prod\limits_{j=1}^{n-1}(3+r_j)+3\left(1-1_{\{i=0\}}\right)\left(\prod\limits_{j\neq i,j=1}^{n-1}(3+r_j)\right)\\
    u_i&\coloneqq p_i-\hat\varepsilon q_i
\end{alignedat}}
\end{equation}
Note that constraint $(C1)$ is actually 2 strict inequality constraints. Next, it is useful to know when the inverse problem has a unique minimizer or if it even does have a unique minimizer. So far, it has been shown that the general problem is quasiconvex. For the problem to have a unique minimizer, strict quasiconvexity is required. It is shown that in general, the inverse problem at hand is not strictly quasiconvex.

First, it makes sense to assess the signs of entries of $\bm{u}$ since it is already known that $\bm\phi$ lives in the simplex ($\mathcal{C}$). From definitions \ref{eqn: main optimization problem full definition}, $\bm{u}=\bm{p}-\hat\varepsilon\bm{q}$. So,

\begin{equation}\label{eqn: U entries}
\begin{array}{rcl}
u_0 &=& \left[\prod\limits_{j=1}^{n-1}(3+r_j)\right](1-\hat\varepsilon)\\[2pt]
u_1 &=& p_{1}-\hat\varepsilon q_{1}
      =3\left(\prod\limits_{j=2}^{n-1}\left(3+r_j\right)\right)\left(1+r_1-\hat\varepsilon\right)
      =3\left(\prod\limits\limits_{j=2}^{n-1}\left(3+r_j\right)\right)\left(\frac{\varepsilon_1}{\varepsilon_0}-\hat\varepsilon\right)\\[2pt]
u_2 &=& p_{2}-\hat\varepsilon q_{2}
      =3\left(\prod\limits_{j=1,j\ne 2}^{n-1}\left(3+r_j\right)\right)\left(\frac{\varepsilon_2}{\varepsilon_0}-\hat\varepsilon\right)\\[2pt]
u_3 &=& p_{3}-\hat\varepsilon q_{3}
      =3\left(\prod\limits_{j=1,j\ne 3}^{n-1}\left(3+r_j\right)\right)\left(\frac{\varepsilon_3}{\varepsilon_0}-\hat\varepsilon\right)\\[4pt]
\multicolumn{3}{c}{\text{\Huge$\vdots$}}\\[4pt]
u_{n-1} &=& p_{n-1}-\hat\varepsilon q_{n-1}=3\left(\prod\limits_{j=1}^{n-2}\left(3+r_j\right)\right)\left(\frac{\varepsilon_{n-1}}{\varepsilon_0}-\hat\varepsilon\right)    
\end{array}
\end{equation}
Also note,
\begin{equation}\label{eqn: ehat in range general n}
    \min\left(1,\frac{\varepsilon_1}{\varepsilon_0},\cdots,\frac{\varepsilon_{n-1}}{\varepsilon_0}\right)<\hat\varepsilon<\max\left(1,\frac{\varepsilon_1}{\varepsilon_0},\cdots,\frac{\varepsilon_{n-1}}{\varepsilon_0}\right);\quad \hat\varepsilon\in\mathrm{Im}\left(\overline\varepsilon\right)
\end{equation}
Now, equations \ref{eqn: U entries} and \ref{eqn: ehat in range general n} show that some entries of $\bm{u}$ are negative and some are positive (based on $\varepsilon_i\lessgtr\hat\varepsilon$). This clearly means that the problem \ref{eqn: main optimization problem full definition} is not strictly quasiconvex since the expression within the absolute value in $\eta$ can cross zero for several $\bm\phi$. It is important to note here that it is unlikely that an entry of $\bm{u}$ will be zero in simple composite materials. For example, in complex composite materials, it is possible that a component is chosen not because of its permittivity properties but rather because of some other desirable property. So, it is important to consider when it is useful to use dielectric constants to study the inverse problem for a given composite material.    

\section{Solutions}\label{section: solutions}

The solution to the optimization problem \ref{eqn: main optimization problem full definition} is $\mathcal{M}$,
\begin{equation}
\begin{gathered}
    \mathcal{M}=\mathcal{C}\>\cap\>\mathcal{H}\text{, with, }\\
    \mathcal{H}\coloneqq\{\bm\phi\>|\>\bm{u}\cdot\bm\phi=0\},\>\mathcal{C}\coloneqq\{\bm\phi\>|\>\bm{1^n}\cdot\bm\phi=1,\>\phi_i>0\>\forall\>i\in\{0,1,2,\cdots,(n-1)\}\}
\end{gathered}
\end{equation}

Instead of the usual simplex ($\mathcal{C}$), it makes more sense to consider $\bm\phi$ to be the dominant component ordered simplex $\left(\mathcal{C}_{(0)}\right)$. Since the matrix material `$0$' will have a volume fraction greater than all other components (that is, $\phi_0\ge\phi_i\>\forall i>0$). For composites with spherical inclusions, the matrix material will have a much higher volume fraction than this, especially if the inclusions are all spherical and of same size (packing fraction, etc). It is important to note here that there is no mention of size of inclusions explicitly anywhere in the Eshelby-Mori-Tanaka model. In fact, this is a drawback of the model (see \cite{BERRYMAN1996149},\cite{Norris1989}, etc. on realizability). So, the actual physically feasible set of $\bm\phi$ will have fewer entries than the ordered simplex. But the results have been shown for the simplex and the ordered simplex in this article and similar calculations can be performed for the specific problem at hand.

To analytically determine the set of all minimizers $\mathcal{M}$, the signs of entries of $\bm{u}$ and the edges of the set in which $\bm\phi$ lives are used. $\eta=0$ on an edge when $\overline\varepsilon-\hat\varepsilon$ goes from negative on one vertex of the set to positive on an adjacent vertex of the set. The signs of the entries of $\bm{u}$ and the edges of the set in which $\bm\phi$ lives are useful in this regard. The set $\mathcal{M}$ is then the convex hull of all all points at which $\eta=0$ on the edges of the set in which $\bm\phi$ lives. This set will automatically be a subset of $\mathcal{C}$. It is possible that a vertex of the set of $\bm\phi$ is itself in $\mathcal{M}$, but this happens only when a particular entry of $\bm{u}$ is zero. In the next 2 subsections, the solution sets $\mathcal{M}$ for $\bm\phi\in\mathcal{C}$ and for $\phi\in\mathcal{C}_{(0)}$ are determined.

\subsection{\texorpdfstring{$\bm{\phi}$ is in the Simplex ($\mathcal{C}$)}{phi in Simplex (C)}}
\begin{equation}
    \mathcal{C}\coloneqq\{\bm\phi\>|\>\bm{1^n}\cdot\bm\phi=1,\>\phi_i>0\>\forall\>i\in\{0,1,2,\cdots,(n-1)\}\}
\end{equation}
The vertices of $\mathcal{C}$ are $\bm{e_i}$ where $\bm{e_i}$ form the standard/Kronecker delta basis ($\bm{e_i}\cdot\bm{e_j}=\delta_{ij}$) for $\mathbb{R}^n$. $\mathcal{C}$ has $n$ vertices and $(n(n-1))/2$ edges. The following definitions are made to state the result clearly:
\begin{equation}
    \mathcal{P}_{1}\coloneqq\{i\>|\>u_i>0\};\>\mathcal{P}_{2}\coloneqq\{i\>|\>u_i<0\};\>\mathcal{Z}\coloneqq\{i\>|\>u_i=0\}
\end{equation}
A vector $\bm\phi^{ij}$ on the edge of $\mathcal{C}$ can be parametrized by $t_{ij}\in(0,1)$,
\begin{equation}\label{eqn: phi edge}
    \bm\phi^{ij}(t_{ij})=t_{ij}\bm{e_i}+(1-t_{ij})\bm{e_j}\>\forall\>i,j\in\{0,1,2,\cdots,(n-1)\},i\neq j
\end{equation}
So, the minimizer set $\mathcal{M}$ becomes:
\begin{equation}
    \mathcal{M}_{\mathcal{C}}=conv\{\{e_i\>|\>i\in\mathcal{Z}\}\>\cup\>\{\bm\phi^{ij\star}\>|\>i\in\mathcal{P}_{1},j\in\mathcal{P}_{2}\}\}\text{ where, }\bm{u}\cdot\bm\phi^{ij\star}=0
\end{equation}
Where $\bm{e_i}$ such that $i\in\mathcal{Z}$ is a minimizer that says the composite material is not a composite material but is not a composite material but just the $i^\text{th}$ material. This solution is ignored. $\bm\phi^{ij\star}$ are $\bm{\phi}$ vectors on the edges of $\mathcal{C}$ such that $\bm{u}\cdot \bm\phi^{ij\star}=0,i\in\mathcal{P}_{1},j\in\mathcal{P}_{2}$. Notice,
\begin{equation}
    \bm{u}\cdot\bm\phi^{ij\star}=t_{ij}^\star u_i+\left(1-t_{ij}^\star \right)u_j=0\text{ so, }t_{ij}^\star =\frac{u_j}{u_j-u_i}
\end{equation}
So, from equation \ref{eqn: phi edge}, one has,
\begin{equation}
    \bm\phi^{ij\star}=\frac{u_j}{u_j-u_i}\bm{e_i}-\frac{u_i}{u_j-u_i}\bm{e_j}
\end{equation}
There will be $|\mathcal{P}_{1}|\cdot|\mathcal{P}_{2}|$ such points (a maximum of $(n(n-1))/2$). Thus, the solution is,
\begin{equation}\label{eqn: phi in simplex}
\boxed{
\mathcal{M}=conv\left\{\bm\phi^{ij\star}\>|\>\bm\phi^{ij\star}=\frac{u_j}{u_j-u_i}\bm{e_i}-\frac{u_i}{u_j-u_i}\bm{e_j}\>\forall\>u_i>0,u_j<0\right\}
}
\end{equation}
$\bm{u}$ entries are computed using equations \ref{eqn: U entries}.

\subsection{\texorpdfstring{$\bm\phi$ is in the Dominant-Component Ordered Simplex $\left(\mathcal{C}_{(0)}\right)$}{phi in the Ordered Simplex (C0)}}
The dominant-component ordered simplex is defined as,
\begin{equation}
    \mathcal{C}_{(0)}=\{\bm\phi\in\mathbb{R}^n:\phi_i\ge0,\bm{1^n}\cdot\bm{\phi}=1,\phi_0\ge\phi_i\>\forall\>i\in\{1,2,3,\cdots,n-1\}\}
\end{equation}
$\mathcal{C}_{(0)}$ has $2^{n-1}$ vertices. $S\subseteq\{1,2,3,\cdots,( n-1)\}$ is used to refer to these vertices:
\begin{equation}
    \phi^S_{i}=\begin{cases}
        \frac{1}{|S|+1}, &i=0 \text{ or }i\in S,\\
        0, &i\notin \{0\}\cup S
    \end{cases}
\end{equation}
As an example, for $n=4$, the $8$ vertices are:
\begin{equation}\label{eqn: n=4 vertices}
\begin{gathered}
    \bm{\phi^{S=\varnothing}}=(1,0,0,0),\quad \bm{\phi^{S=\{1\}}}=\left(\frac{1}{2},\frac{1}{2},0,0\right),\quad  
    \bm{\phi^{S=\{2\}}}=\left(\frac{1}{2},0,\frac{1}{2},0\right),\\
    \bm{\phi^{S=\{3\}}}=\left(\frac{1}{2},0,0,\frac{1}{2}\right),\quad \bm{\phi^{S=\{1,2\}}}=\left(\frac{1}{3},\frac{1}{3},\frac{1}{3},0\right),\quad \bm{\phi^{S=\{2,3\}}}=\left(\frac{1}{3},0,\frac{1}{3},\frac{1}{3}\right)\\
    \bm{\phi^{S=\{1,3\}}}=\left(\frac{1}{3},\frac{1}{3},0,\frac{1}{3}\right),\quad\bm{\phi^{S=\{1,2,3\}}}=\left(\frac{1}{4},\frac{1}{4},\frac{1}{4},\frac{1}{4}\right)
\end{gathered}
\end{equation}
$\mathcal{C}_{(0)}$ has $(n-1)2^{n-2}$ edges. So, for $n=4$, there are $12$ edges. These are:
\begin{equation}
    \begin{gathered}
        \bm{\phi^{S=\varnothing}}-\bm{\phi^{S=\{1\}}},\quad \bm{\phi^{S=\varnothing}}-\bm{\phi^{S=\{2\}}},\quad \bm{\phi^{S=\varnothing}}-\bm{\phi^{S=\{3\}}},\\
        \bm{\phi^{S=\{1\}}}-\bm{\phi^{S=\{1,2\}}},\quad \bm{\phi^{S=\{1\}}}-\bm{\phi^{S=\{1,3\}}},\\\bm{\phi^{S=\{2\}}}-\bm{\phi^{S=\{1,2\}}},\quad \bm{\phi^{S=\{2\}}}-\bm{\phi^{S=\{2,3\}}},\\
        \bm{\phi^{S=\{3\}}}-\bm{\phi^{S=\{1,3\}}},\quad \bm{\phi^{S=\{3\}}}-\bm{\phi^{S=\{2,3\}}},\\
        \bm{\phi^{S=\{1,2\}}}-\bm{\phi^{S=\{1,2,3\}}},\quad \bm{\phi^{S=\{2,3\}}}-\bm{\phi^{S=\{1,2,3\}}},\quad
        \bm{\phi^{S=\{1,3\}}}-\bm{\phi^{S=\{1,2,3\}}}
    \end{gathered}
\end{equation}
For the case of $\bm\phi\in\mathcal{C}$, the vertex values of $\bm{u}\cdot\bm\phi$ were the entries of $\bm{u}$ themselves. This is not the case for $\bm\phi\in\mathcal{C}_{(0)}$. Denote $u^S$ to be $\bm{u}\cdot\bm\phi^S$ at the vertex associated with $S\subseteq\{1,2,3,\cdots,(n-1)\}$. Then, from equations \ref{eqn: U entries},
\begin{equation}\label{eqn: u^S value}
    u^S\coloneqq \bm{u}\cdot\bm{\phi^S};\quad u^S=\frac{u_0+\sum\limits_{i\in S}u_i}{|S|+1}
\end{equation}
$u^S$ and $u^{S\cup\{k\}}$ with $k\notin S$ are $\bm{u}\cdot\bm{\phi^S}$ values at adjacent vertices. $\bm{\phi^{S,k}}$ is used to denote $\bm \phi$ values on the edge joining the vertices $\bm{\phi^S}$ and $\bm{\phi^{S\cup\{k\}}}$. So,
\begin{equation}
\begin{gathered}
    \bm{\phi^{S,k}}=t_{S,k}\bm{\phi^{S}}+(1-t_{S,k})\bm{\phi^{S\cup\{k\}}};\\
    t_{S,k}\in(0,1),k\in\{1,2,\cdots,(n-1)\},k\notin S,S\subseteq\{1,2,\cdots,(n-1)\}
\end{gathered}
\end{equation}
If the above equation for $\bm{\phi^{S,k}}$ is dotted with $\bm{u}$ and set to zero, one has,
\begin{equation}
    \bm{u}\cdot\bm{\phi^{S,k\star}}=t_{S,k}^\star u^S+(1-t_{S,k}^\star )u^{S\cup\{k\}}=0;\quad t_{S,k}^\star =\frac{u^{S\cup\{k\}}}{u^{S\cup\{k\}}-u^S}
\end{equation}
Note, that in this case, only the vertex $\bm\phi=(1/(|S|+1),1/(|S|+1),\cdots,1/(|S|+1))$ is a composite material which comprises all components. Other vertices have at least one zero entry. So $\bm{u}\cdot\bm\phi$ being zero at a vertex is ignored and not stated explicitly. Either way, the solution written as a convex hull of edge crossings $\bm{\phi^{S,k\star}}$ will contain this value of $\bm\phi$ if it is a solution. So,
\begin{equation}
\boxed{
\begin{gathered}
    \mathcal{M}_{\mathcal{C}_{(0)}}=conv\left\{\bm{\phi^{S,k\star}}\left\lvert\bm{\phi^{S,k\star}}=t^\star _{S,k}\bm{\phi^S}+\left(1-t^\star _{S,k}\right)\bm{\phi^{S\cup\{k\}}}\>\forall\>u^S u^{S\cup\{k\}}<0\right.\right\}\\
S\subseteq\{1,2,\cdots,(n-1)\},\>k\in\{1,2,\cdots,(n-1)\},k\notin S\\
\phi^S_{i}=\begin{cases}
        \frac{1}{|S|+1}, &i=0 \text{ or }i\in S,\\
        0, &i\notin \{0\}\cup S
    \end{cases};\quad u^S=\frac{u_0+\sum\limits_{i\in S}u_i}{|S|+1};\quad t_{S,k}^\star =\frac{u^{S\cup\{k\}}}{u^{S\cup\{k\}}-u^S}
\end{gathered}
}
\end{equation}

\subsection{\texorpdfstring{Unique Minimizer $\bm\phi^\star $}{Unique Minimizer phi*}}
So, it is shown that the inverse Eshelby-Mori-Tanaka problem for isotropic composites with spherical inclusions reduces to determining a single $\bm\phi$ in a convex polytope $\mathcal{M}$. In convex optimization, the focus is on finding numerical solutions of utility rather than analytic closed form solutions. So, in what follows, an optimization problem is shown which can be solved using a solver like cvxpy. Features of the optimization problem will decide on the quality and usefulness of the solution computed using it. These will be discussed in fair detail. There are many ways that this can be done, but it must follow the physics in question. First, an elementary method is shown in this subsection as an example that will not perform well.

\begin{equation}\label{eqn: Ordered Simplex Barrier Problem}
    \begin{aligned}
        \underset{\bm{\phi}}{\text{Minimize: }}-\sum\limits_i\ln(\phi_i)-\sum\limits_{i\ge1}\ln(\phi_0-\phi_i)\quad\text{subject to:}\quad \bm{1}^n\cdot\bm{\phi}=1,\>\bm{u}\cdot\bm{\phi}=0
    \end{aligned}
\end{equation}

Problem \ref{eqn: Ordered Simplex Barrier Problem} above is not a good formulation to solve for the minimizer $\bm{\phi^\star }$ since it is a logarithmic centroid of the set of minimizers $\mathcal{M}_{\mathcal{C}_{(0)}}$ that depends only on the geometry of the dominant-component ordered simplex $\mathcal{C}_{(0)}$ and not any of the physics incorporating vectors like $\bm{u},\bm{p},$ or $\bm{q}$. The second constraint in the problem is just a halfplane and does not help in determining a good physics abiding minimizer from the polytope $\mathcal{M}_{\mathcal{C}_{(0)}}$.

\section{Optimization Problems with the Physics Encoded in the Constraints}\label{section: Final 2 formulations}
\subsection{Epigraph Formulation with Charnes Cooper Change of Variables}\label{subsection: PD-IPM formulation}


\begin{equation}\label{eqn: epigraph with primal-dual interior pt method inclusions only}
\boxed{
\begin{alignedat}{3}
\underset{\bm{z},s,t}{\text{Minimize: }}&t\\
\text{Subject to: } &\bm{y}\cdot\bm{z}+p_0s=1 &&(C1)\\
&\underbrace{\bm{x}\cdot\bm{z}+p_0s}_{=\overline{\varepsilon}}\le t+\hat\varepsilon &&(C2)\\
-&\underbrace{\left(\bm{x}\cdot\bm{z}+p_0s\right)}_{=\overline\varepsilon}\le t-\hat\varepsilon &&(C3)\\
&\bm{1^{(n-1)}}\cdot\bm{z}\le s; &&(C4) \\
&\left(s-\bm{1^{(n-1)}}\cdot\bm{z}\right)\ge z_k\>\forall\>k\in\{1,2,\cdots,(n-1)\} &&(C5)\\
&\bm{z}\succeq \bm{0};\quad s\ge 0,\quad t\ge0; && (C6)-(C8)\\
r_i\coloneqq&\frac{\varepsilon_i}{\varepsilon_0}-1\>\forall i\in[1,\cdots,(n-1)];\quad p_0\coloneqq\prod\limits_{k=1}^{n-1}(3+r_k)\\
    p_i^{inc}\coloneqq& 3\left((1+r_i)\prod\limits_{j\neq i,j=1}^{n-1}(3+r_j)\right);\quad \bm{x}\coloneqq\bm{p^{inc}}-p_0\bm{1^{(n-1)}};\\    q_i^{inc}\coloneqq&3\left(\prod\limits_{j\neq i,j=1}^{n-1}(3+r_j)\right);\quad \bm{y}\coloneqq \bm{q^{inc}}-p_0\bm{1^{(n-1)}};\\
    \bm{\phi^{inc}}\coloneqq& \left(\phi_1,\phi_2,\cdots,\phi_{n-1}\right);\quad \bm{z}\coloneqq s\bm{\phi^{inc}};\quad s\coloneqq\frac{1}{\bm{y}\cdot\bm{\phi^{inc}}+p_0};
\end{alignedat}}
\end{equation}
This problem is a Linear Programming (LP) problem that is arrived at using the epigraph formulation and the Charnes-Cooper change of variables (see Section 4.3.2. in \cite{boyd2004convex}). An LP is known to be convex, so a proof is not required for convexity.  

However, though Problem \ref{eqn: epigraph with primal-dual interior pt method inclusions only} has correct physics enforcing constraints, it is still not enough since there are $(n-1)$ independent variables (components of $\bm{\phi}$) to determine and not enough independent constraints. Constraint $(C1)$ of Problem \ref{eqn: epigraph with primal-dual interior pt method inclusions only} enforces change of variables which is to do with algebra and $(C2),$ $(C3)$ are essentially just one constraint from the physics. At least $(n-1)$ \textit{independent constraints} that enforce the physics are required to determine $\bm{\phi^\star}$. However, a problem formulation using the Charnes-Cooper change of variables such as problem \ref{eqn: epigraph with primal-dual interior pt method inclusions only} is not ideal to incorporate multi-frequency data, as the scaling constraint enforcing the change of variables $(C1)$ causes problems. Thus, another formulation is developed in the following section.

\subsection{Multi-Frequency, Complex Dielectric Constant and Epigraph based Formulation}\label{subsection: Final Formulation}

One way to deal with the lack of independent constraints is to measure the dielectric constant (real part) at at least $(n-1)$ different frequencies for all dispersive phases. One can also measure complex dielectric constant (real and imaginary parts) and obtain independent constraints. Problem \ref{eqn: main optimization problem full definition} must be modified to reflect this change. In problem \ref{eqn: main optimization problem multifrequency full definition}, subscript (superscript) ``$(k)$" is used for scalars (vectors) to signify measurement at different frequencies.  

\begin{equation}\label{eqn: main optimization problem multifrequency full definition}
\boxed{
\begin{alignedat}{3}
\underset{\bm{\phi},\eta}{\text{Minimize: }}&\eta\coloneqq\underset{k}{\mathrm{max}}\left\vert\overline{\varepsilon}_{(k)}-\hat{\varepsilon}_{(k)}\right\vert=\underset{k}{\mathrm{max}}\left\lvert\frac{\bm{p^{(k)}}\cdot\bm{\phi}}{\bm{q^{(k)}}\cdot\bm{\phi}}-\hat{\varepsilon}_{(k)}\right\rvert\\
\text{Subject to: }\\
&\bm{\phi}\succeq\bm{0};\quad\bm{1^n}\cdot\bm{\phi}=1& \\
\hline
    1_{\{i=0\}}&\coloneqq\begin{cases}
        0,i>0,\\
        1,i=0
    \end{cases};\> r_{i(k)}\coloneqq\frac{\varepsilon_{i(k)}}{\varepsilon_0}-1\>\forall i\in\{1,\cdots,(n-1)\};\\
    k&\in\{1,\cdots,m\}\text{ with }m\ge(n-1); \text{ (measurement frequencies) }\\
    \hline
    p_i^{(k)}\coloneqq \left(1_{\{i=0\}}\right)\prod\limits_{i=1}^{n-1}&(3+r_{i(k)})+3\left(1-1_{\{i=0\}}\right)\left(\left(1+r_{i(k)}\right)\prod\limits_{j\neq i,j=1}^{n-1}\left(3+r_{i(k)}\right)\right);\\
    q_i^{(k)}\coloneqq\left(1_{\{i=0\}}\right)\prod\limits_{i=1}^{n-1}&\left(3+r_{i(k)}\right)+3\left(1-1_{\{i=0\}}\right)\left(\prod\limits_{j\neq i,j=1}^{n-1}\left(3+r_{i(k)}\right)\right)
\end{alignedat}}
\end{equation}

The epigraph of $\eta$ is then used. This will now give more independent constraints. The following robust one-shot LP using the epigraph formulation of problem \ref{eqn: main optimization problem multifrequency full definition}, considering all components except the matrix material to be dispersive, and considering complex dielectric constants is formulated:

\begin{equation}\label{eqn: Final formulation - Multifreq Epigraph}
\boxed{
\begin{alignedat}{3}
\underset{\bm{\phi},t,\bm{c},\bm{d}}{\text{Minimize: }}&t\\
\text{Subject to: } 
&\bm{\phi}\succeq\bm{0};\quad\bm{1^n}\cdot\bm{\phi}=1 ;\quad t\ge0;\quad\quad \quad &&(C1)-(C3)\\
-&\bm{A}\bm{\phi}\preceq \bm{c};\quad \bm{A}\bm{\phi}\preceq \bm{c}   &&(C4)\\
-&\bm{B}\bm{\phi}\preceq \bm{d};\quad \bm{B}\bm{\phi}\preceq \bm{d} &&(C5)\\
&\bm{c}+\bm{d}\preceq t\bm{q^{\mathrm{max}}} &&(C6)\\
\hline
&k\in\{1,\cdots,m\}(m\text{ frequencies})&&\\
&\bm{a^{(k)}}\coloneqq\mathfrak{R}\left(\bm{p^{(k)}-\hat\varepsilon_{(k)}q^{(k)}}\right);\quad &&\bm{b^{(k)}}\coloneqq\mathfrak{I}\left(\bm{p^{(k)}-\hat\varepsilon_{(k)}q^{(k)}}\right);\\
&\bm{A}\coloneqq\begin{pmatrix}
    \bm{a^{(1)}}&\bm{a^{(2)}}&\cdots&\bm{a^{(m)}}
\end{pmatrix}^\top;\quad &&\bm{B}\coloneqq\begin{pmatrix}
    \bm{b^{(1)}}&\bm{b^{(2)}}&\cdots&\bm{b^{(m)}}
\end{pmatrix}^\top; \\
&\bm{q^{\mathrm{max}}}\coloneqq\underset{i}{\mathrm{max}}\left\lvert q^{(k)}_i\right\rvert; \quad\bm{A},\bm{B}\in\mathbb{R}^{m\times n}; &&\quad\bm{c},\bm{d},\bm{q^\mathrm{max}}\in\mathbb{R}^{m} 
\end{alignedat}}
\end{equation}
Here, in problem \ref{eqn: Final formulation - Multifreq Epigraph}, $(C4)$ and $(C5)$ are $m$ independent constraints each that enforce the physics. Note that the component by which the dielectric constants are normalized (for $r_{i(k)}$ for each component) must be non-dispersive (no change in dielectric constant with frequency of measurement). So in problem \ref{eqn: Final formulation - Multifreq Epigraph}, $\varepsilon_0$ must be non-dispersive. If not, another component that is non-dispersive must be used to normalize. If all the components of the composite material are dispersive, then a straightforward reformulation of the problem can be done. This choice was made just to simplify the presentation here. Notice that this formulation is capable of taking dielectric constant measurements of the composite material at multiple frequencies unlike problem \ref{eqn: epigraph with primal-dual interior pt method inclusions only}. 

How would one know exactly what number and kind of constraints $(C4)$ and $(C5)$ are \textit{sufficient} to uniquely and correctly identify $\bm{\phi^\star}$? For this, a \textit{constraint sensitivity matrix} is determined for the system in Problem \ref{eqn: Final formulation - Multifreq Epigraph}. A simplex tangent subspace is first defined, since for an $n$ component composite material, the number of independent variables are $n-1$.
\begin{equation}
    \begin{gathered}
        \bm{\phi^1},\bm{\phi^2}\in\mathcal{C},\text{ then }\bm{1^n}\cdot\bm{\Delta\phi}=0\text{ where }\bm{\Delta\phi}=\left(\bm{\phi^1}-\bm{\phi^2}\right). \\
        \text{ So define tangent subspace: }\mathcal{T}\coloneqq\{\bm{\Delta\phi}\in\mathbb{R}^n:\bm{1^n}\cdot\bm{\Delta\phi}=0\}
    \end{gathered}
\end{equation}
Now, define an orthonormal basis $\bm{V}$ that spans $\mathcal{T}$:  
\begin{equation}\label{eqn: basis spanning T}
    \begin{gathered}
        \bm{S}\coloneqq\begin{pmatrix}
            \bm{I^{(n-1)}} & -\bm{1^{(n-1)}}
        \end{pmatrix}^\top\text{ and }\bm{V}\coloneqq\bm{S}\left(\bm{S}^\top\bm{S}\right)^{-\frac{1}{2}};\>\bm{S},\bm{V}\in\mathbb{R}^{n\times(n-1)}\\
        \text{Notice, }\bm{\Delta\phi}=\bm{V\Delta\phi^{inc}} \text{ and } \left\lVert\bm{\Delta\phi}\right\rVert_2=\left\lVert\bm{\Delta\phi^{inc}}\right\rVert_2
    \end{gathered}
\end{equation}
Thus, the constraint sensitivity matrix $\bm{G}$ that operates on $\bm{\Delta\phi^{inc}}$ can be definied as follows:
\begin{equation}\label{eqn: G matrix}
    \begin{gathered}
        \bm{G}\coloneqq\begin{pmatrix}
            \bm{QA}\\\bm{QB}
        \end{pmatrix}\bm{V};\quad\bm{Q}\coloneqq\mathrm{diag}\left(1/q^\mathrm{max}_1, 1/q^\mathrm{max}_2,\cdots,1/q^\mathrm{max}_{m} \right);\quad\bm{G}\in\mathbb{R}^{2m\times(n-1)};\\
        \text{So, }\begin{pmatrix}
            \bm{QA}\\\bm{QB}
        \end{pmatrix}\bm{\Delta\phi}=\bm{G\Delta\phi^{inc}}
    \end{gathered}
\end{equation}
So, in order for the inverse to be identifiable, one must have $\mathrm{Rank}(\bm{G})$ at least $n-1$ for the identifiability of $\bm{\phi^\star}$. This means that if only real parts of the dielectric constants are measured, then measurements at at least $m=n-1$ frequencies are required. If complex dielectric constants are measured, then measurements at at least $m=(n-1)/2$ frequencies are required. Another measure that determines the quality of predicted volumetric split is the lowest singular value of $\bm{G}$ ($\sigma_\mathrm{min}(\bm{G})$). In general, the higher the value of $\sigma_\mathrm{min}(\bm{G})$, the better the quality of the split. This is shown by deriving a crude bound on $\left\lVert\bm{\phi^\star}-\bm{\phi^\mathrm{t}}\right\rVert_\infty$ after considering the measurements to be noisy.

In Problem \ref{eqn: Final formulation - Multifreq Epigraph}, $\hat \varepsilon_{(k)}$ are assumed to be noisy as would be the case in real experiments. So, $\hat{\varepsilon}_{(k)}$ for all $k$ are such that:
\begin{equation}
\begin{gathered}
\hat\varepsilon_{(k)} \coloneqq\hat\varepsilon_{\mathrm{t}(k)}+\hat\varepsilon_{\mathrm{n}(k)};\quad \hat\varepsilon_{\mathrm{n}(k)} = \varepsilon_{\mathrm{n}(k)}/\varepsilon_0;\quad\varepsilon_{\mathrm{n}(k)}=\varepsilon_{\mathrm{n}(k)_{R}}+i\varepsilon_{\mathrm{n}(k)_{I}};\\
\varepsilon_{\mathrm{n}(k)_{R}}\sim\mathcal{U}\left(-\delta_R,\delta_R\right);\quad \varepsilon_{\mathrm{n}(k)_{I}}\sim\mathcal{U}\left(-\delta_I,\delta_I\right)    
\end{gathered}
\end{equation}
Here, $\hat\varepsilon_{\mathrm{t}(k)}$ refers to the true normalized dielectric constant of the composite material measured at frequency $k$ (total $m$ frequencies). $\hat\varepsilon_{\mathrm{n}(k)}$ is the noise component of $\hat\varepsilon_{(k)}$. $\left(\bm{\phi^\star},t^\star,\bm{c^\star},\bm{d^\star}\right)$ is the optimizer of problem \ref{eqn: Final formulation - Multifreq Epigraph} and $\bm{\phi^\mathrm{t}}$ is the true volumetric split of the composite material associated with the dielectric constants $\hat\varepsilon_{\mathrm{t}(k)}$ at $m$ different frequencies. $\bm{\Delta\phi}$ is defined to be $\bm{\phi^\star}-\bm{\phi^\mathrm{t}}$. Some new definitions are made to determine bounds on noise given a desired $\left\lVert\bm{\phi}\right\rVert_2$.

\begin{gather}
    \bm{u^{\hat{\varepsilon}_{\mathrm{t}}(k)}}\coloneqq\bm{p^{(k)}}-\hat{\varepsilon}_{\mathrm{t}(k)}\bm{q^{(k)}};\quad\text{Notice, }\bm{u^{\hat{\varepsilon}_{\mathrm{t}}(k)}}\cdot\bm{\phi^\mathrm{t}}=0.\label{eqn: u true dot phi true is zero}\\
    \bm{u^{\hat\varepsilon(k)}}\coloneqq\bm{p^{(k)}}-\hat{\varepsilon}_{(k)}\bm{q^{(k)}}= \bm{u^{\hat{\varepsilon}_{\mathrm{t}}(k)}}-\hat\varepsilon_{\mathrm{n}(k)}\bm{q^{(k)}}  \label{eqn: u true and noise definitions};\\
    \bm{u^{\hat\varepsilon_{\mathrm{t}}(k)}}\cdot\bm{\Delta\phi}=\bm{u^{\hat\varepsilon_{\mathrm{t}}(k)}}\cdot\bm{\phi^\star}-\underbrace{\bm{u^{\hat\varepsilon_{\mathrm{t}}(k)}}\cdot\bm{\phi^\mathrm{t}}}_{=0}=\bm{u^{\hat\varepsilon_{\mathrm{t}}(k)}}\cdot\bm{\phi^\star}\label{eqn: u delta phi equals u phi star};\\
    \text{Using }\bm{a^{(k)}},\bm{b^{(k)}},\bm{A},\text{ and }\bm{B}\text{ from Problem \ref{eqn: Final formulation - Multifreq Epigraph} }:\nonumber\\
    \bm{a^{\hat{\varepsilon}_{\mathrm{t}}(k)}}\coloneqq\mathfrak{R}\left(\bm{u^{\hat{\varepsilon}_{\mathrm{t}}(k)}}\right);\>\bm{b^{\hat{\varepsilon}_{\mathrm{t}}(k)}}\coloneqq\mathfrak{I}\left(\bm{u^{\hat{\varepsilon}_{\mathrm{t}}(k)}}\right)\text{ and so on.}\\
    \bm{A^{\hat{\varepsilon}_{\mathrm{t}}}}\coloneqq\begin{pmatrix}
        \bm{a^{\hat{\varepsilon}_{\mathrm{t}}(1)}}&\cdots&\bm{a^{\hat{\varepsilon}_{\mathrm{t}}(m)}}
    \end{pmatrix}^\top;\quad     \bm{B^{\hat{\varepsilon}_{\mathrm{t}}}}\coloneqq\begin{pmatrix}
        \bm{b^{\hat{\varepsilon}_{\mathrm{t}}(1)}}&\cdots&\bm{b^{\hat{\varepsilon}_{\mathrm{t}}(m)}}
    \end{pmatrix}^\top
\end{gather}
The following identities are used:
\begin{gather}
    |z_R+iz_I|\le|z_R|+|z_I|;\> z_R,z_I\in\mathbb{R};\tag{I1}\label{eqn: complex I1}\\
    \left\lvert\mathfrak{R}(z)\right\rvert+\left\lvert\mathfrak{I}(z)\right\rvert\le2|z|;\>z,w\in\mathbb{C};\tag{I2}\label{eqn: complex I2}\\ 
    \left\lvert\mathfrak{R}(z+w)\right\rvert+\left\lvert\mathfrak{I}(z+w)\right\rvert\le \left(\left\lvert\mathfrak{R}(z)\right\rvert+\left\lvert\mathfrak{I}(z)\right\rvert\right)+\left(\left\lvert\mathfrak{R}(w)\right\rvert+\left\lvert\mathfrak{I}(w)\right\rvert\right)\tag{I3}\label{eqn: complex I3}
\end{gather}    
Using identity \ref{eqn: complex I1} and the definition of $\hat\varepsilon_{\mathrm{n}(k)}$, one has
\begin{equation}
    \hat\varepsilon_{\mathrm{n}}^{\mathrm{max}}\coloneqq\underset{k}{\mathrm{max}}\left\lvert\hat{\varepsilon}_{\mathrm{n(k)}}\right\rvert\le\frac{\delta_R+\delta_I}{\varepsilon_0}
\end{equation}
Here, $t^\star$ is the minimized objective. Next, dotting equation \ref{eqn: u true and noise definitions} with $\bm{\phi^\star}$, and using \ref{eqn: u delta phi equals u phi star},
\begin{equation}
    \bm{u^{\hat\varepsilon_{\mathrm{t}}(k)}}\cdot\bm{\Delta\phi}=\bm{u^{\hat\varepsilon_{\mathrm{t}}(k)}}\cdot\bm{\phi^\star}=\bm{u^{\hat\varepsilon(k)}}\cdot\bm{\phi^\star}+\hat\varepsilon_{\mathrm{n}(k)}\left(\bm{q^{(k)}}\cdot\bm{\phi^\star}\right)
\end{equation}
Next use identity \ref{eqn: complex I3} to get:
\begin{equation}
\begin{aligned}
    \left\lvert\bm{a^{\hat\varepsilon_{\mathrm{t}}(k)}}\cdot\bm{\Delta\phi}\right\rvert+\left\lvert\bm{b^{\hat\varepsilon_{\mathrm{t}}(k)}}\cdot\bm{\Delta\phi}\right\rvert&\le
    \left[\left\lvert\bm{a^{\hat\varepsilon(k)}}\cdot\bm{\phi^\star}\right\rvert +\left\lvert\bm{b^{\hat\varepsilon(k)}}\cdot\bm{\phi^\star}\right\rvert \right]\\
    &+\left[\left\lvert\mathfrak{R}\left( \hat\varepsilon_{\mathrm{n(k)}}\bm{q^{(k)}}\cdot\bm{\phi^\star}\right)\right\rvert +\left\lvert\mathfrak{I}\left( \hat\varepsilon_{\mathrm{n(k)}}\bm{q^{(k)}}\cdot\bm{\phi^\star}\right)\right\rvert \right]\label{eqn: utrue dot phistar le}
\end{aligned}
\end{equation}
For a bound on the first $\left[\cdot\right]$, constraint $(C6)$ of problem \ref{eqn: Final formulation - Multifreq Epigraph} is used. For a bound on the second $\left[\cdot\right]$, identity \ref{eqn: complex I2} is used.
\begin{gather}
    \left\lvert\bm{a^{\hat\varepsilon(k)}}\cdot\bm{\phi^\star}\right\rvert+\left\lvert\bm{b^{\hat\varepsilon(k)}}\cdot\bm{\phi^\star}\right\rvert\le t^\star q^{\mathrm{max}}_k;\label{eqn: u ehat dot phi star bound}\\
    \left\lvert\mathfrak{R}\left( \hat\varepsilon_{\mathrm{n(k)}}\bm{q^{(k)}}\cdot\bm{\phi^\star}\right)\right\rvert +\left\lvert\mathfrak{I}\left( \hat\varepsilon_{\mathrm{n(k)}}\bm{q^{(k)}}\cdot\bm{\phi^\star}\right)\right\rvert\le2\hat\varepsilon_{\mathrm{n}}^{\mathrm{max}} q^{\mathrm{max}}_k \label{eqn: noise dot phistar le}
\end{gather}
Putting inequalities \ref{eqn: utrue dot phistar le}, \ref{eqn: u ehat dot phi star bound}, and \ref{eqn: noise dot phistar le}, one has for each $k$,
\begin{equation}
    \left\lvert\bm{a^{\hat\varepsilon_{\mathrm{t}}(k)}}\cdot\bm{\Delta\phi}\right\rvert+\left\lvert\bm{b^{\hat\varepsilon_{\mathrm{t}}(k)}}\cdot\bm{\Delta\phi}\right\rvert\le\left(t^\star+2\hat\varepsilon_{\mathrm{n}}^\mathrm{max}\right)q^{\mathrm{max}}_k
\end{equation}
Divinding through $q^\mathrm{max}_k$, stacking all inequalities up using the constraint sensitivity matrix \ref{eqn: G matrix}, and using the definition of the $L_2$ norm,
\begin{equation}
\begin{aligned}
    \left\lVert\begin{pmatrix}
    \bm{QA^{\hat\varepsilon_{\mathrm{t}}}} \\ \bm{QB^{\hat\varepsilon_{\mathrm{t}}}}
\end{pmatrix}\bm{\Delta\phi}\right\rVert_2=\left\lVert\bm{G^{\hat\varepsilon_{\mathrm{t}}}\Delta\phi^{inc}}\right\rVert_2=\sigma_\mathrm{min}\left(\bm{G^{\hat\varepsilon_{\mathrm{t}}}}\right)\left\lVert\bm{\Delta\phi}\right\rVert_2&\le\sqrt{2m}\left(t^\star+2\hat\varepsilon_{\mathrm{n}}^\mathrm{max}\right)\\
&\le \sqrt{2m}\left[t^\star+2\left(\frac{\delta_R+\delta_I}{\varepsilon_0}\right)\right]
\end{aligned}
\end{equation}
Here, $\sigma_\mathrm{min}(\bm{G})$ is the minimum singular value of $\bm{G}$. Apart from $\mathrm{Rank}(\bm{G})$, the other main value to look at is $\sigma_\mathrm{min}(\bm{G})$. This determines the quality of the split and depends on the material system. Given observed noise $\left(\delta_R+\delta_I\right)$, the quality of the determined volumetric split may be determined:
\begin{equation}\label{eqn: Bound phi norm error}
\boxed{
\left\lVert\bm{\Delta\phi}\right\rVert_\infty\le\left\lVert\bm{\Delta\phi}\right\rVert_2\le\frac{ \sqrt{2m}}{\sigma_\mathrm{min}\left(\bm{G^{\hat\varepsilon_{\mathrm{t}}}}\right)}\left[t^\star+2\left(\frac{\delta_R+\delta_I}{\varepsilon_0}\right)\right]}
\end{equation}
The purpose of arriving at the above bound is to show that the frequencies at which measurements are to be made can be decided by examining the matrix $\bm{G^{\hat\varepsilon_t}}$ and $t^\star$ values for a system. The requirement for $t^\star$ to be low is an obvious requirement. Apart from this, the lowest singular value of $\bm{G}$, $\sigma_\mathrm{min}\left(\bm{G^{\hat\varepsilon_t}}\right)$, must be a high value. This value depends on $\bm{\phi^\mathrm{t}}$, but for a material system, if any one component is very dispersive in a given frequency range and the measurements are taken in this range, then $\sigma_\mathrm{min}\left(\bm{G^{\hat\varepsilon_t}}\right)$ will be a high value for that system. This bound is not a tight bound and depends on $\bm{\phi^\mathrm{t}}$ and is not helpful beyond showing the importance of the constraint sensitivity matrix $\bm{G}$. It shows that $\bm{G^{\hat\varepsilon_t}}$ and $\bm{t^\star}$ are important to examine to ensure the constraints employed are sufficient to obtain a good volumetric split.
In the next section, Problem \ref{eqn: Final formulation - Multifreq Epigraph} will be applied to 3 different 3-component material systems, and the behavior of this inequality, and through it, the behavior of Problem \ref{eqn: Final formulation - Multifreq Epigraph}, will be studied. In this paper, only complex dielectric constant measurements are considered. If other physical information is known about the system, $\bm{G}$ can be modified and its rank and $\sigma_\mathrm{min}$ can be determined to see if problem \ref{eqn: Final formulation - Multifreq Epigraph} will work well. 

In this work, for all computations, cvxpy with a Primal-Dual Interior Point Method (P-D IPM) based solver called Embedded Conic Optimization Solver (ECOS, see \cite{ECOS}) is used. A brief overview of the working of P-D IPMs for LPs is given in chapter 11 of \cite{boyd2004convex}. ECOS is a fast solver built for embedded applications and can potentially be used for real-time volumetric splits of composites.

\section{Applications}\label{section: application}

In this section, Problem \ref{eqn: Final formulation - Multifreq Epigraph} is applied to three different 3 component material systems. Each material system has only one dispersive component and 2 non-dispersive components. 
\begin{itemize}
    \item Material System 1 (MS-1): Epoxy, Glass Microspheres, Pores.
    \item Material System 2 (MS-2): Aggregates, Portland Cement Paste, Pores.
    \item Material System 3 (MS-3): Carbon-Loaded Epoxy, Glass Microspheres, Pores.
\end{itemize}
Ideally, the $\hat{\hat{\varepsilon}}_{\mathrm{t}(k)}$ values are supposed to be experimentally determined. However, to validate the performance of the inverse method, the $\hat{\hat{\varepsilon}}_{\mathrm{t}(k)}$ values in this work are computed using the forward model \ref{eqn: MTB} and an appropriate dispersion model for the dispersive component in each material system. Then, noise as described in the previous section is added to obtain $\hat\varepsilon_{(k)}$ which is then passed to the optimization problem \ref{eqn: Final formulation - Multifreq Epigraph}. Several randomly generated $\bm{\phi^\mathrm{t}}$ values are used in the validation. The plots comparing volume fractions (in Figures \ref{fig: MS-1, MS-3} and \ref{fig: MS-2}) consider 100 randomly generated $\bm{\phi^\mathrm{t}}$ values. The plots comparing errors across the number of frequencies at which measurements are taken ($m$) consider 30,000 randomly generated $\bm{\phi^\mathrm{t}}$ values. The randomly generated volumetric splits ($\bm{\phi^\mathrm{t}}$) are such that realizability is non-issue. For simplicity, a noise level of $\delta_R=\delta_I=0.1$ is considered in all computations in this section.   

The frequency band for the ``measured" values is considered on the basis of the dispersion properties of the dispersive component. For MS-1 and MS-2 the 0.4-3.0 GHz range is chosen. This is also the range of commercially available Ground Penetrating Radar (GPR) equipment. For MS-3, the frequency band 0.5-1.0 MHz was used. 

For the first and third material systems listed, the microstructure is almost exactly the one considered in this paper. The glass microspheres are manufactured to be exactly spheres and the pores are assumed to be spherical. For the Portland Cement Concrete system, an equivalent microstructure is required. This will be discussed in the next section.

\subsection{MS-1: Epoxy, Glass, Pores}

The microstructure details are as follows:
\begin{itemize}
    \item Neat Epoxy material is the matrix material subscripted by `$0$'. A commonly used model (single Debye relaxation model with conductivity (see \cite{debye1929polar}) for its dispersion in the 0.4-3.0 GHz band is
    \begin{equation}
    \begin{gathered}
        \varepsilon_0(\omega)=\varepsilon_{0(\infty)}+\frac{\Delta\varepsilon_0}{1-i\omega \tau_0}+\frac{is_0}{\omega\varepsilon_{\mathrm{vacuum}}}\text{ with,}\\
        \varepsilon_{0(\infty)} = 2.90,\> \Delta\varepsilon_0 = 0.6,\> \tau_0 = 1.5\times 10^{-10}\mathrm{s},\> s_{0} = 10^{-12},\\\varepsilon_{\mathrm{vacuum}}=8.854187817\times10^{-12}\mathrm{F/m}
    \end{gathered}
    \end{equation}
    For $m=1$, the frequency 2.0 GHz is considered. For $m=2,3,4,5$, equally spaced frequencies in the band 0.4-3.0 GHz are considered. 
    \item Glass appears as spherical inclusions in the microstructure and is subscripted by `$1$'. The dielectric constant of glass is taken to be a constant $\varepsilon_1=5.5 + 0.05i$.
    \item Air voids or pores appear in spheres and air is subscripted by `$2$'. The dielectric constant of air is taken to be a constant $\varepsilon_2=1.0006$.
    \item $\bm{\phi^\mathrm{t}}\in\{\bm{\phi^\mathrm{t}}|\>\phi_0\ge0.65,\>\phi_1\ge0.1, \phi_2\ge0.02,\>\bm{1^3}\cdot\bm{\phi^\mathrm{t}}=1\}$
\end{itemize}
The above properties are used with model \ref{eqn: MTB} to determine the true dielectric constant at $m$ different frequencies. Noise is then added as modeled in subsection \ref{subsection: Final Formulation} before passing it to Problem \ref{eqn: Final formulation - Multifreq Epigraph} to determine the volumetric split. Essentially, the goal is to determine $\bm{\phi^\star}$ such that it is equal to $\bm{\phi^\mathrm{t}}$. This procedure is followed for the other two material systems studied in this section as well. 

Neat epoxy is not very dispersive and this is seen clearly as a consequence in the effectiveness of Problem \ref{eqn: Final formulation - Multifreq Epigraph} in predicting the volumetric split in Figure \ref{fig: MS-1, MS-3}. Increasing the number of frequencies of measurement $m$ does not seem to improve performance of Problem \ref{eqn: Final formulation - Multifreq Epigraph} as illustrated in Figures \ref{fig: phi norm delta MS-1}, \ref{fig: tstar delta MS-1}, \ref{fig: min sigma min G MS-1}.

\begin{figure}[!htbp]
\centering

\setcounter{subfigure}{0}
\renewcommand{\thesubfigure}{a(\arabic{subfigure})}

\begin{subfigure}[t]{0.45\textwidth}
  \centering
  \includegraphics[width=\linewidth]{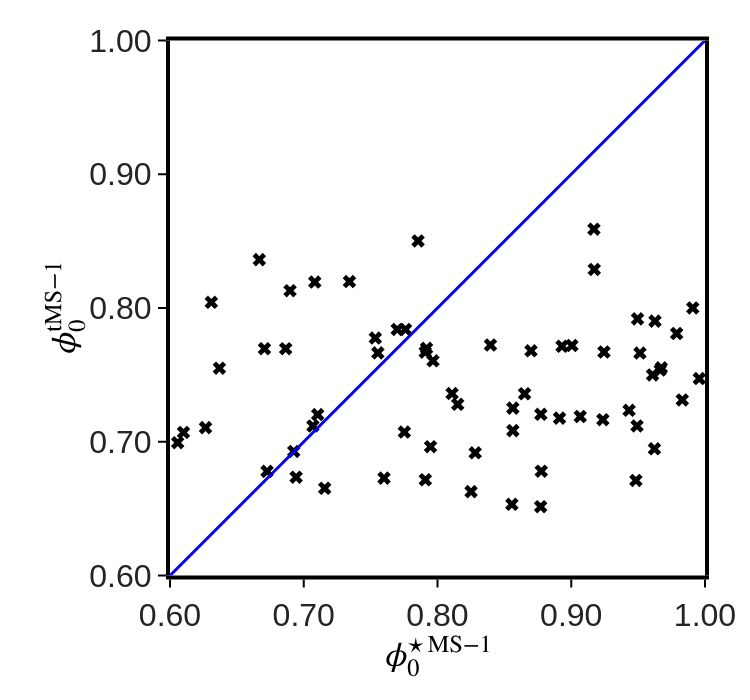}
  \caption{Volume Fraction of neat Epoxy in MS-1 for $\delta_R=\delta_I=0.1$ and for $m=1$ (measurement at 2.0 GHz).}
  \label{fig: EGP-Epoxy}
\end{subfigure}\hfill%
\begin{subfigure}[t]{0.45\textwidth}
  \centering
  \includegraphics[width=\linewidth]{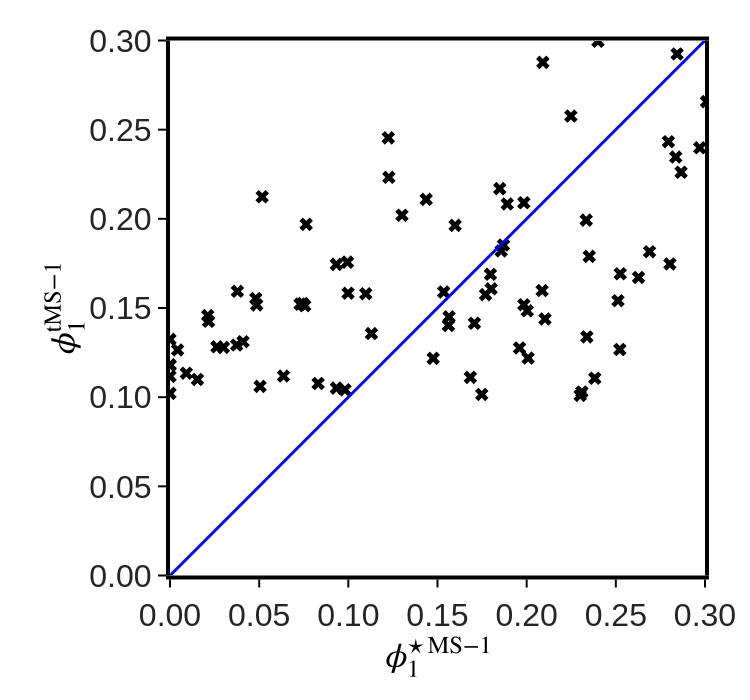}
  \caption{Volume Fraction of Glass in MS-1 for $\delta_R=\delta_I=0.1$ and for $m=1$ (measurement at 2.0 GHz).}
  \label{fig: EGP-Glass}
\end{subfigure}

\vspace{4pt}

\begin{subfigure}[t]{0.45\textwidth}
  \centering
  \includegraphics[width=\linewidth]{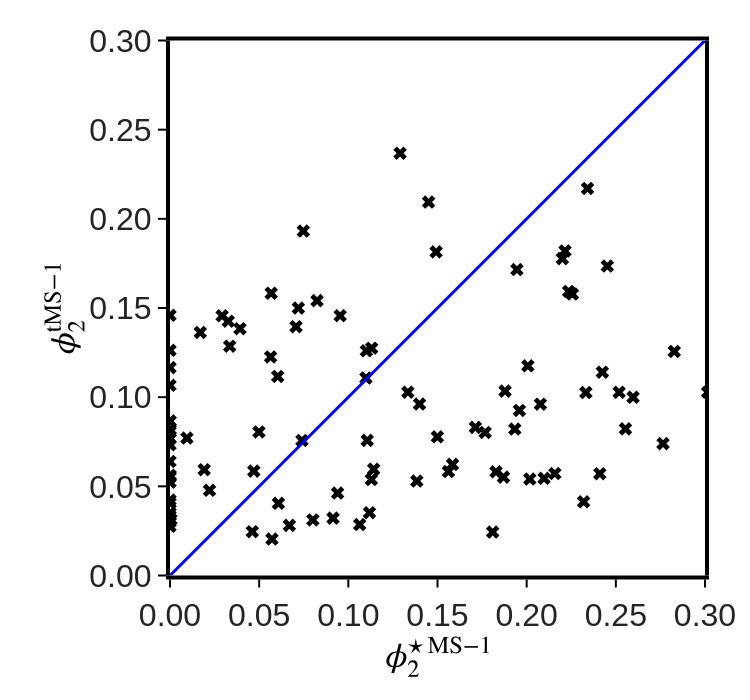}
  \caption{Volume Fraction of Pores in MS-1 for $\delta_R=\delta_I=0.1$ and for $m=1$ (measurement at 2.0 GHz).}
  \label{fig: EGP-Pores}
\end{subfigure}\hfill%
\setcounter{subfigure}{0}
\renewcommand{\thesubfigure}{b\arabic{subfigure}}
\begin{subfigure}[t]{0.45\textwidth}
  \centering
  \includegraphics[width=\linewidth]{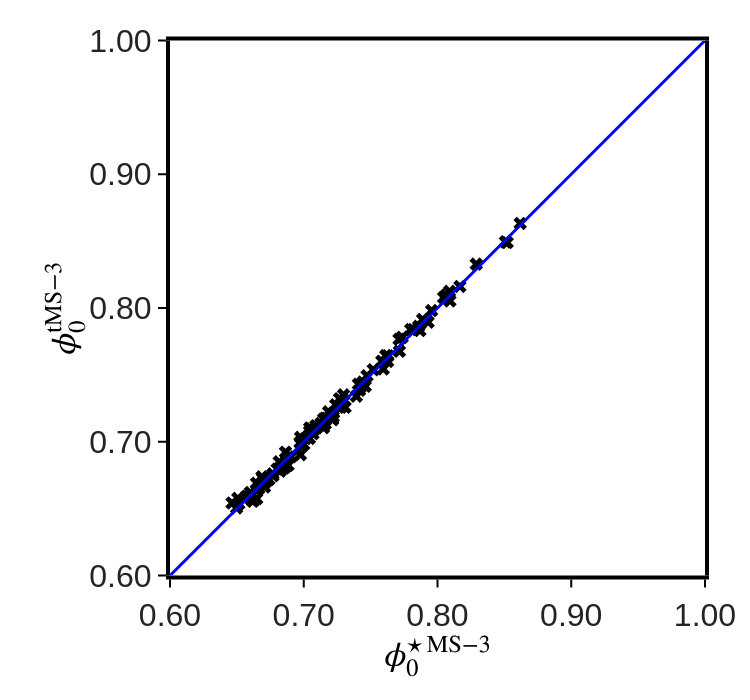}
  \caption{Volume Fraction of Carbon-Loaded Epoxy in MS-3 for $\delta_R=\delta_I=0.1$ and for $m=1$ (measurement at 0.5 MHz).}
  \label{fig: CLEGP-CLE}
\end{subfigure}

\vspace{4pt}

\begin{subfigure}[t]{0.45\textwidth}
  \centering
  \includegraphics[width=\linewidth]{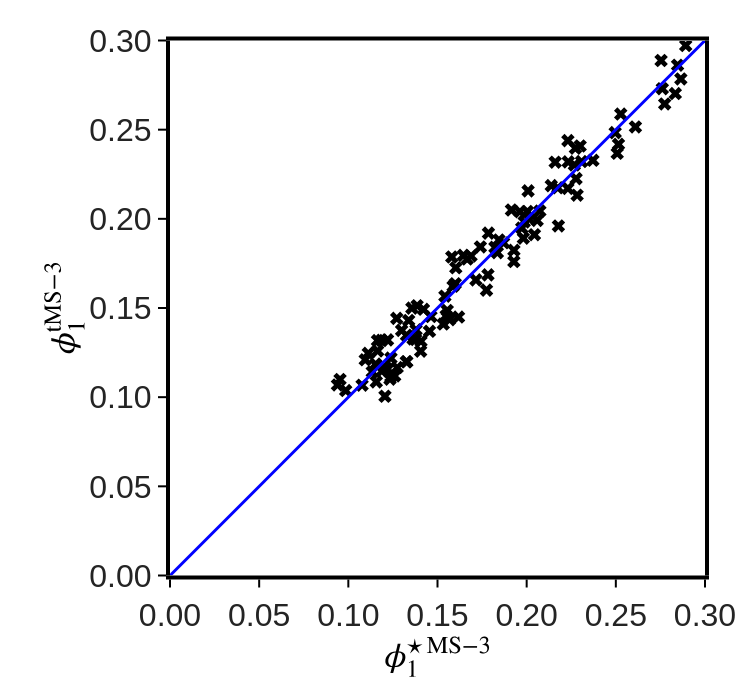}
  \caption{Volume Fraction of Glass in MS-3 for $\delta_R=\delta_I=0.1$ and for $m=1$ (measurement at 0.5 MHz).}
  \label{fig: CLEGP-Glass}
\end{subfigure}\hfill%
\begin{subfigure}[t]{0.45\textwidth}
  \centering
  \includegraphics[width=\linewidth]{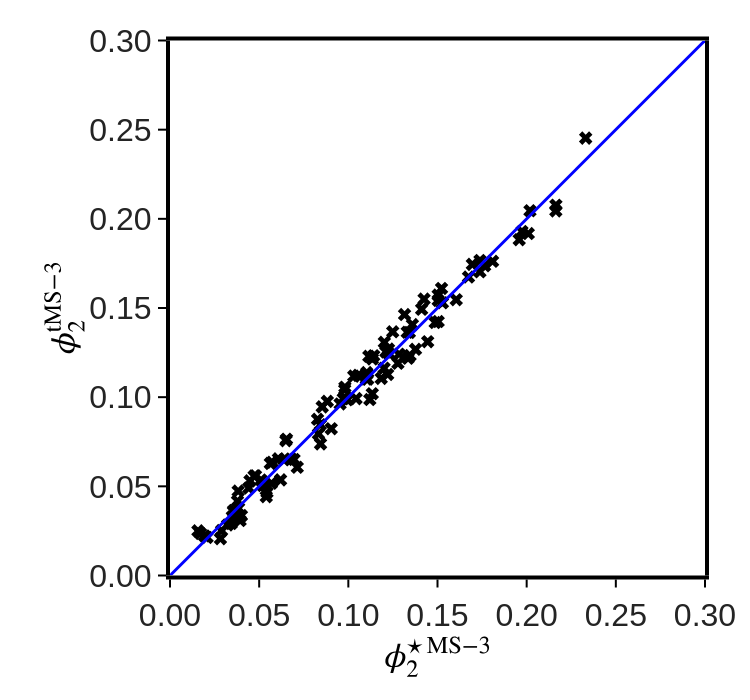}
  \caption{Volume Fraction of Pores in MS-3 for $\delta_R=\delta_I=0.1$ and for $m=1$ (measurement at 0.5 MHz).}
  \label{fig: CLEGP-Pores}
\end{subfigure}

\caption{Performance of Problem \ref{eqn: Final formulation - Multifreq Epigraph} in terms of recovery of $\bm{\phi^\mathrm{t}}$ for MS-1 and MS-3 for $\delta_R=\delta_I=0.1$ and for $m=1$. For MS-1, measurement frequency is $2.0$ GHz and for MS-3, measurement frequency is $0.5$ MHz.}
\label{fig: MS-1, MS-3}
\end{figure}

\subsection{MS-2: Aggregates, Portland Cement Paste, Pores}
The actual microstructure for this material system is very complex. Ideally, to use an equivalent microstructure as simple as the one that is described in what follows, different material properties ($\varepsilon_0^\star, \varepsilon_1^\star, \varepsilon_2^\star$) than the actual ($\varepsilon_0, \varepsilon_1, \varepsilon_2$) must be used (see next section). However, this subsection assumes that ($\varepsilon_0^\star, \varepsilon_1^\star, \varepsilon_2^\star$) are the same as ($\varepsilon_0, \varepsilon_1, \varepsilon_2$) to show how the system behaves. The details of the equivalent microstructure are as follows:
\begin{itemize}
    \item Aggregate is the matrix material subscripted by `$0$'. The dielectric constant of aggregate is taken to be a constant $\varepsilon_1=5.5 + 0.05i$.
    \item Portland cement paste appears as spherical inclusions in the microstructure and is subscripted by `$1$'. This component is almost as dispersive as neat epoxy with respect to the real part of dielectric constant but much more dispersive with respect to the imaginary part. The Portland cement paste is assumed to be room-conditioned with some moisture in the pores and the parameters are chosen accordingly. A commonly used model (Cole-Cole relaxation model with DC conductivity, see \cite{cole1941dispersion}) for its dispersion in the 0.4-3.0 GHz band is:
    \begin{equation}
    \begin{gathered}
        \varepsilon_1(\omega)=\varepsilon_{1(\infty)}+\frac{\Delta\varepsilon_1}{1+(-i\omega \tau_1)^{1-\beta_1}}+\frac{is_{1(dc)}}{\omega\varepsilon_{\mathrm{vacuum}}}\text{ with,}\\
        \varepsilon_{1(\infty)} = 4.2,\>\Delta\varepsilon_1 = 0.8,\>\tau_1 = 2.0\times10^{-10}\mathrm{s},\>\beta_1 = 0.35,\>s_{1(dc)} = 2.0\times10^{-3}\mathrm{S/m}.
    \end{gathered}
    \end{equation}
    For $m=1$, the frequency 2.0 GHz is considered. For $m=2,3,4,5$, equally spaced frequencies in the band 0.4-3.0 GHz are considered. 
    \item Air voids or pores appear in spheres and air is subscripted by `$2$'. The dielectric constant of air is taken to be a constant $\varepsilon_2=1.0006$.
    \item $\bm{\phi^\mathrm{t}}\in\{\bm{\phi^\mathrm{t}}|\>\phi_0\ge0.65,\>\phi_1\ge0.1, \phi_2\ge0.02,\>\bm{1^3}\cdot\bm{\phi^\mathrm{t}}=1\}$
\end{itemize}

\begin{figure}[!htbp]
\setcounter{subfigure}{0}
\renewcommand{\thesubfigure}{a(\arabic{subfigure})}
    \centering

    \begin{subfigure}[t]{0.45\textwidth}
        \centering
        \includegraphics[width=\linewidth]{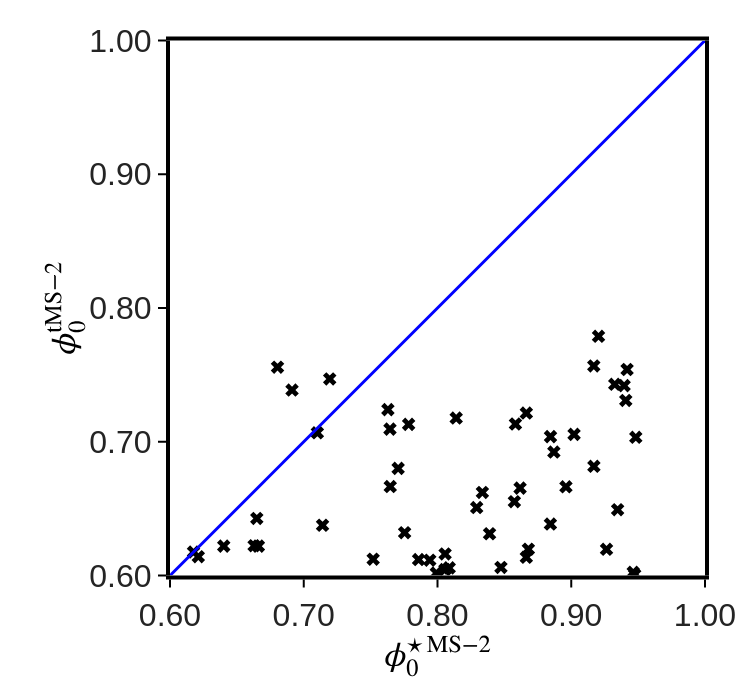}
        \caption{Volume Fraction of Aggregates in MS-2 for $\delta_R=\delta_I=0.1$ and for $m=1$ (measurement at 2.0 GHz).}
        \label{fig: PCC-Agregates m=1}
    \end{subfigure}
    \hfill
    \begin{subfigure}[t]{0.45\textwidth}
        \centering
        \includegraphics[width=\linewidth]{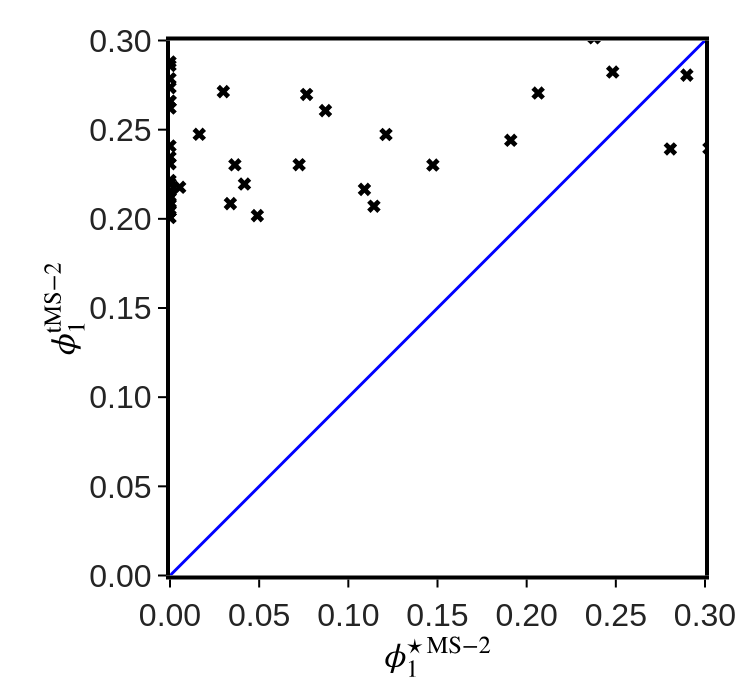}
        \caption{Volume Fraction of Portland Cement Paste in MS-2 for $\delta_R=\delta_I=0.1$ and for $m=1$ (measurement at 2.0 GHz).}
        \label{fig: PCC-PC m=1}
    \end{subfigure}

    \vspace{4pt} 

    \begin{subfigure}[t]{0.45\textwidth}
        \centering
        \includegraphics[width=\linewidth]{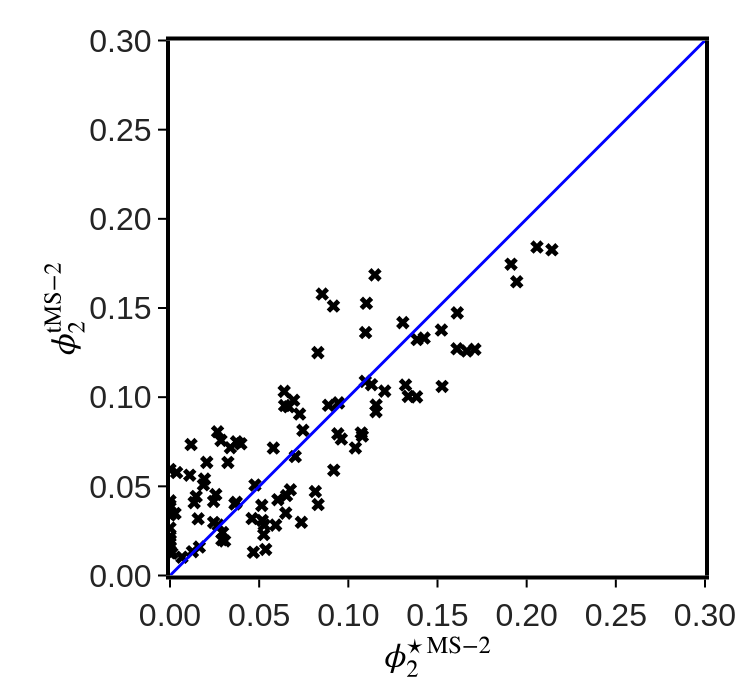}
        \caption{Volume Fraction of Pores in MS-2 for $\delta_R=\delta_I=0.1$ and for $m=1$ (measurement at 2.0 GHz).}
        \label{fig: PCC-Pores m=1}
    \end{subfigure}\hfill
    \setcounter{subfigure}{0}
\renewcommand{\thesubfigure}{b(\arabic{subfigure})}
    \begin{subfigure}[t]{0.45\textwidth}
        \centering
        \includegraphics[width=\linewidth]{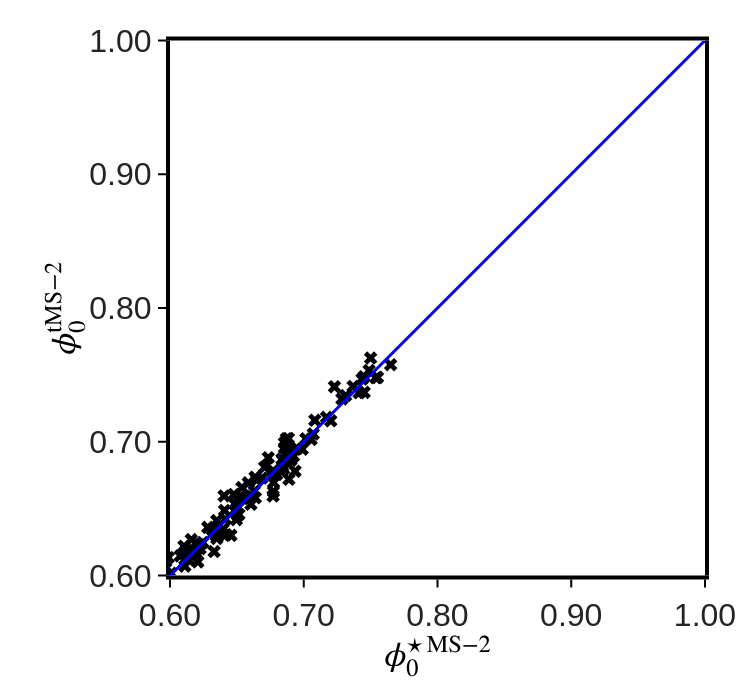}
        \caption{Volume Fraction of Aggregates in MS-2 for $\delta_R=\delta_I=0.1$ and for $m=5$ (measurements at 5 equally spaced frequencies in the band 0.4-3.0 GHz, including 0.4 GHz and 3.0 GHz).}
        \label{fig: PCC-Aggregates m=5}
    \end{subfigure}

    \vspace{4pt}

    \begin{subfigure}[t]{0.45\textwidth}
        \centering
        \includegraphics[width=\linewidth]{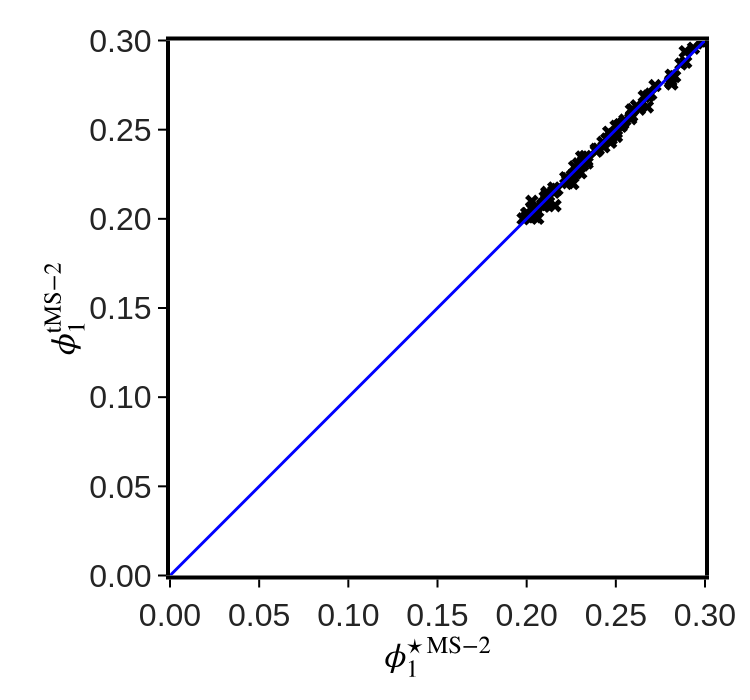}
        \caption{Volume Fraction of Portland Cement Paste in MS-2 for $\delta_R=\delta_I=0.1$ and for $m=5$ (measurements at 5 equally spaced frequencies in the band 0.4-3.0 GHz, including 0.4 GHz and 3.0 GHz).}
        \label{fig: PCC-PC m=5}
    \end{subfigure}
    \hfill
    \begin{subfigure}[t]{0.45\textwidth}
        \centering
        \includegraphics[width=\linewidth]{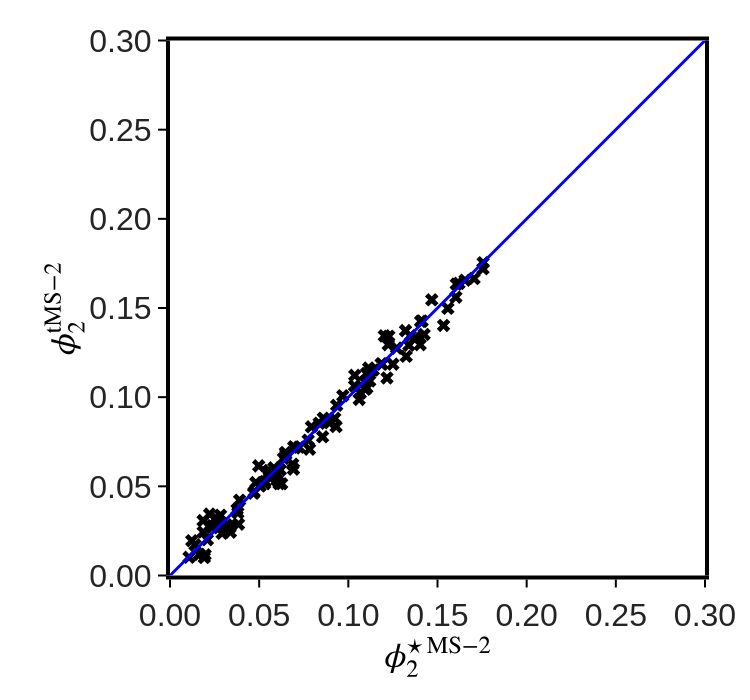}
        \caption{Volume Fraction of Pores in MS-2 for $\delta_R=\delta_I=0.1$ and for $m=5$ (measurements at 5 equally spaced frequencies in the band 0.4-3.0 GHz, including 0.4 GHz and 3.0 GHz).}
        \label{fig: PCC-Pores m=5}
    \end{subfigure}

    \caption{Performance of Problem \ref{eqn: Final formulation - Multifreq Epigraph} (in terms of recovery of $\bm{\phi^\mathrm{t}}$) for MS-2 for $\delta_R=\delta_I=0.1$ and for $m=1$ (measurement at only a frequency of 2.0 GHz) and $m=5$ (measurements at 5 equally spaced frequencies in the band 0.4-3.0 GHz, including 0.4 GHz and 3.0 GHz).}
    \label{fig: MS-2}
\end{figure}

Room-conditioned Portland cement paste with some moisture in the pores is mildly dispersive with respect to the real part of dielectric constant and moderately dispersive with respect to the imaginary part of dielectric constant. In Figure \ref{fig: MS-2}, it is seen that for a noise level of $\delta_R=\delta_I=0.1$, Problem \ref{eqn: Final formulation - Multifreq Epigraph} performs well for $m=5$ and badly for $m=1$ (see Figures \ref{fig: MS-2}, \ref{fig: phi norm delta MS-2}, \ref{fig: tstar delta MS-2}, and \ref{fig: min sigma min G MS-2}). However, notice that the maximum $t^\star$ value is high for $m=5$. In fact, for all the material systems considered, higher $m$ means higher maximum $t^\star$. Also note that for all the material systems considered, higher $m$ also means higher minimum $\sigma_\mathrm{min}\left(\bm{G^{\hat\varepsilon_t}}\right)$.

\subsection{MS-3: Carbon-Loaded Epoxy, Glass Microspheres, Pores}
\begin{figure}[!htbp]
    \centering
\setcounter{subfigure}{0}
\renewcommand{\thesubfigure}{a(\arabic{subfigure})}
    \begin{subfigure}[t]{0.42\textwidth}
        \centering
        \includegraphics[width=\linewidth]{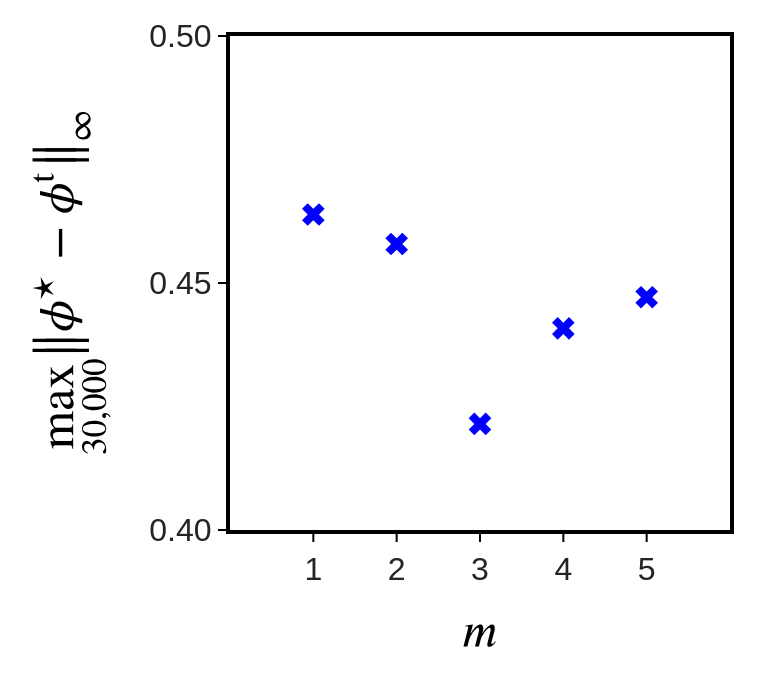}
        \caption{Maximum error in $\bm{\phi}$ for $\delta_R=\delta_I=0.1$ for MS-1. $m=1$: one measurement at 2 GHz. $m=2,3,4,5$: $m$ measurements at equally spaced frequencies in the band 0.4-3.0 GHz, including 0.4 and 3.0 GHz.}
        \label{fig: phi norm delta MS-1}
    \end{subfigure}\hfill
    \begin{subfigure}[t]{0.42\textwidth}
        \centering
        \includegraphics[width=\linewidth]{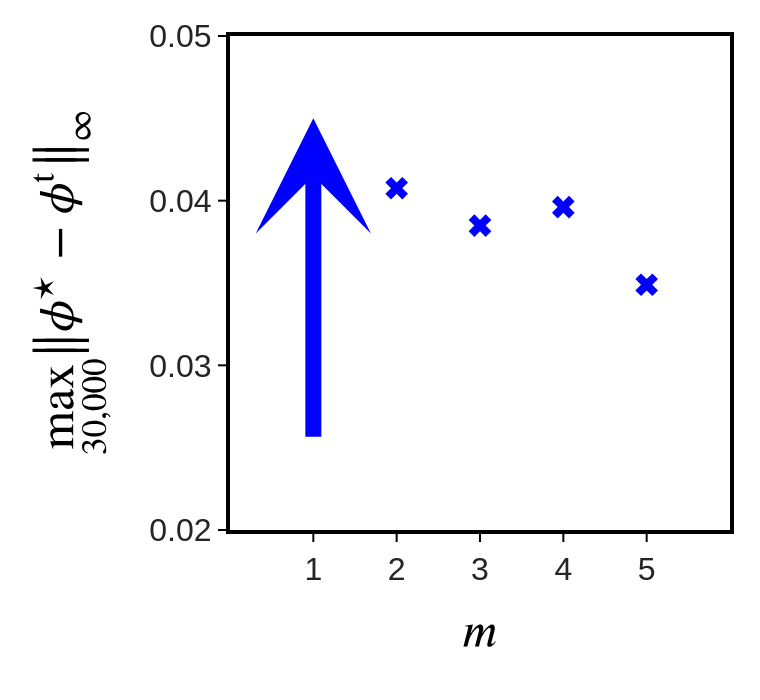}
        \caption{Maximum error in $\bm{\phi}$ for $\delta_R=\delta_I=0.1$ for MS-2. $m=1$: one measurement at 2 GHz. $m=2,3,4,5$: $m$ measurements at equally spaced frequencies in the band 0.4-3.0 GHz, including 0.4 and 3.0 GHz. For $m=1$, $\underset{30,000}{\mathrm{max}}\left\lVert\bm{\phi^\star}-\bm{\phi^\mathrm{t}}\right\rVert_\infty$ is very high $\approx0.62025$.} 
        \label{fig: phi norm delta MS-2}
    \end{subfigure}

    \vspace{4pt} 

    \begin{subfigure}[t]{0.42\textwidth}
        \centering
        \includegraphics[width=\linewidth]{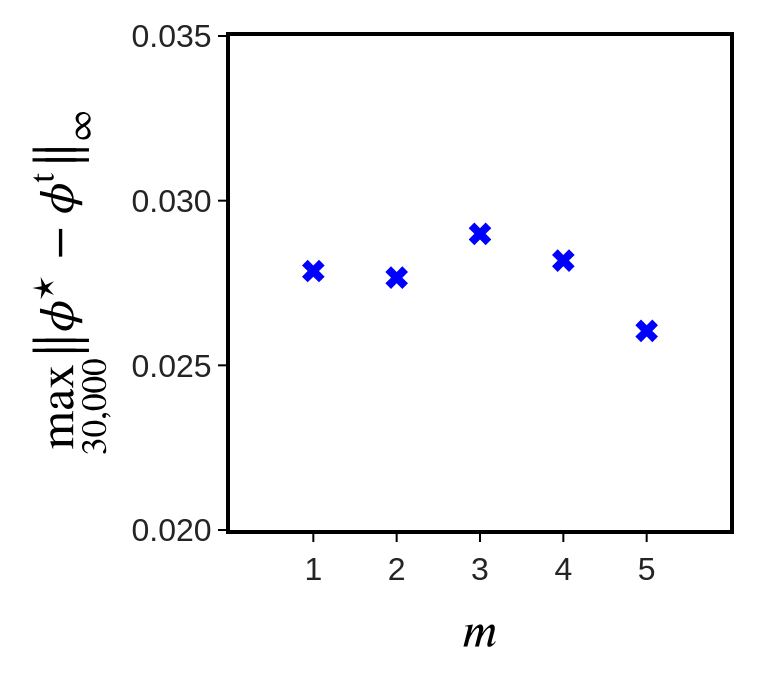}
        \caption{Maximum error in $\bm{\phi}$ for $\delta_R=\delta_I=0.1$ for MS-3. $m=1$: one measurement at 0.5 MHz. $m=2,3,4,5$: $m$ measurements at equally spaced frequencies in the band 0.5-1.0 MHz, including 0.5 and 1.0 MHz.}
        \label{fig: phi norm delta MS-3}
    \end{subfigure}\hfill
    \setcounter{subfigure}{0}
\renewcommand{\thesubfigure}{b(\arabic{subfigure})}
    \begin{subfigure}[t]{0.42\textwidth}
        \centering
        \includegraphics[width=\linewidth]{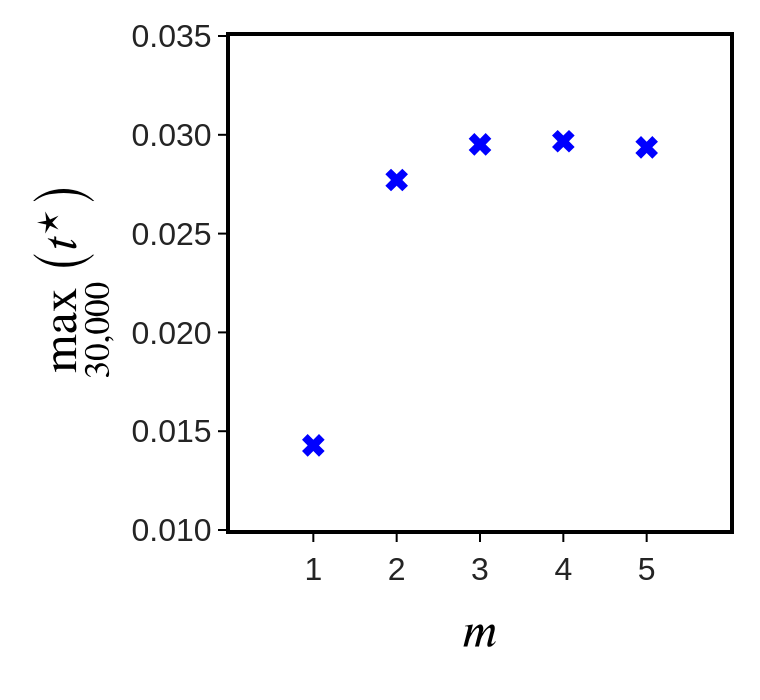}
        \caption{Maximum $t^\star$ for $\delta_R=\delta_I=0.1$ for MS-1. $m=1$: one measurement at 2 GHz. $m=2,3,4,5$: $m$ measurements at equally spaced frequencies in the band 0.4-3.0 GHz, including 0.4 and 3.0 GHz.}
        \label{fig: tstar delta MS-1}
    \end{subfigure}

    \vspace{4pt}

    \begin{subfigure}[t]{0.42\textwidth}
        \centering
        \includegraphics[width=\linewidth]{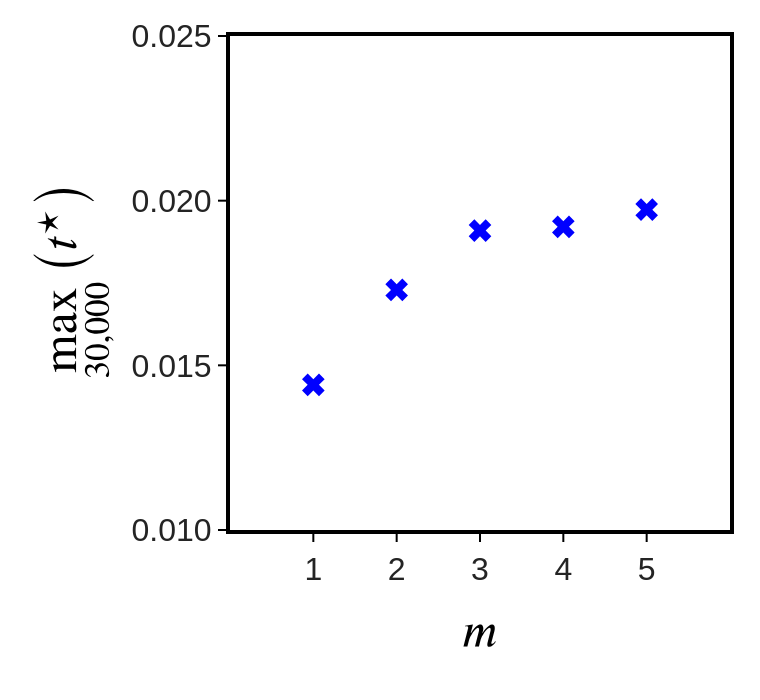}
        \caption{Maximum $t^\star$ for $\delta_R=\delta_I=0.1$ for MS-2. $m=1$: one measurement at 2 GHz. $m=2,3,4,5$: $m$ measurements at equally spaced frequencies in the band 0.4-3.0 GHz, including 0.4 and 3.0 GHz.}
        \label{fig: tstar delta MS-2}
    \end{subfigure}
    \hfill
    \begin{subfigure}[t]{0.42\textwidth}
        \centering
        \includegraphics[width=\linewidth]{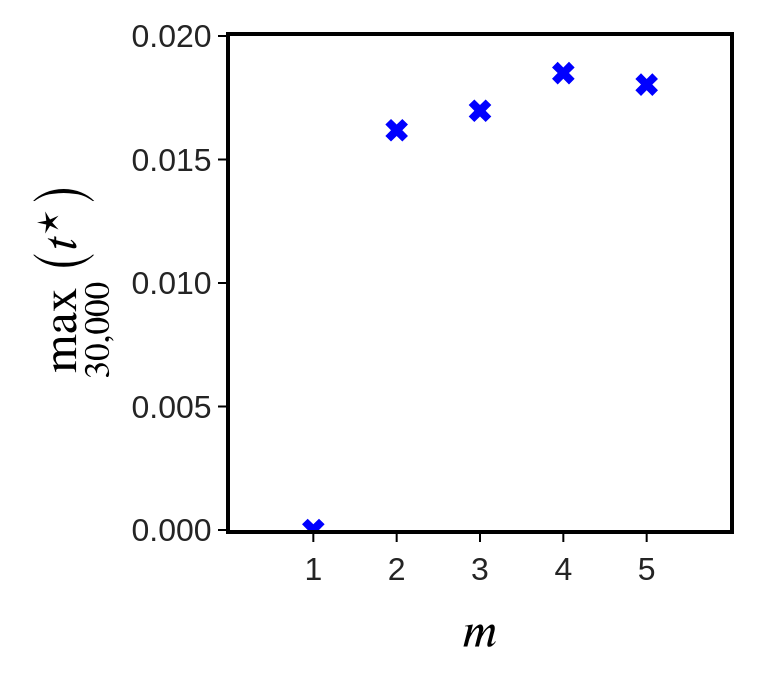}
        \caption{Maximum $t^\star$ for $\delta_R=\delta_I=0.1$ for MS-3. $m=1$: one measurement at 0.5 MHz. $m=2,3,4,5$: $m$ measurements at equally spaced frequencies in the band 0.5-1.0 MHz, including 0.5 and 1.0 MHz.}
        \label{fig: tstar delta MS-3}
    \end{subfigure}

    \caption{Performance of Problem \ref{eqn: Final formulation - Multifreq Epigraph} in terms of $\underset{30,000}{\mathrm{max}}\left\lVert\bm{\phi^\star}-\bm{\phi^\mathrm{t}}\right\rVert_\infty$ and $\underset{30,000}{\mathrm{max}}\left(t^\star\right)$ against number of frequencies at which measurements are taken ($m$) for $\delta_R=\delta_I=0.1$.}
    \label{fig: Problem Performance phi norm error and tstar}
\end{figure}

\begin{figure}[!htbp]
    \centering
\setcounter{subfigure}{0}
\renewcommand{\thesubfigure}{a(\arabic{subfigure})}
    \begin{subfigure}[t]{0.42\textwidth}
        \centering
        \includegraphics[width=\linewidth]{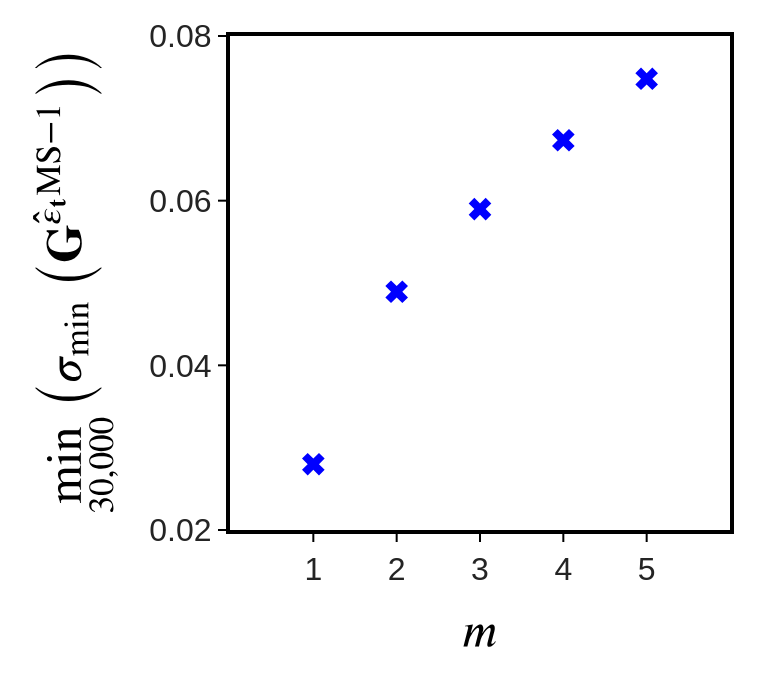}
        \caption{Minimum $\sigma_\mathrm{min}\left(\bm{G^{\hat\varepsilon_t}}\right)$ for $\delta_R=\delta_I=0.1$ for MS-1. $m=1$: one measurement at 2 GHz. $m=2,3,4,5$: $m$ measurements at equally spaced frequencies in the band 0.4-3.0 GHz, including 0.4 and 3.0 GHz.}
        \label{fig: min sigma min G MS-1}
    \end{subfigure}
    \hfill
    \begin{subfigure}[t]{0.42\textwidth}
        \centering
        \includegraphics[width=\linewidth]{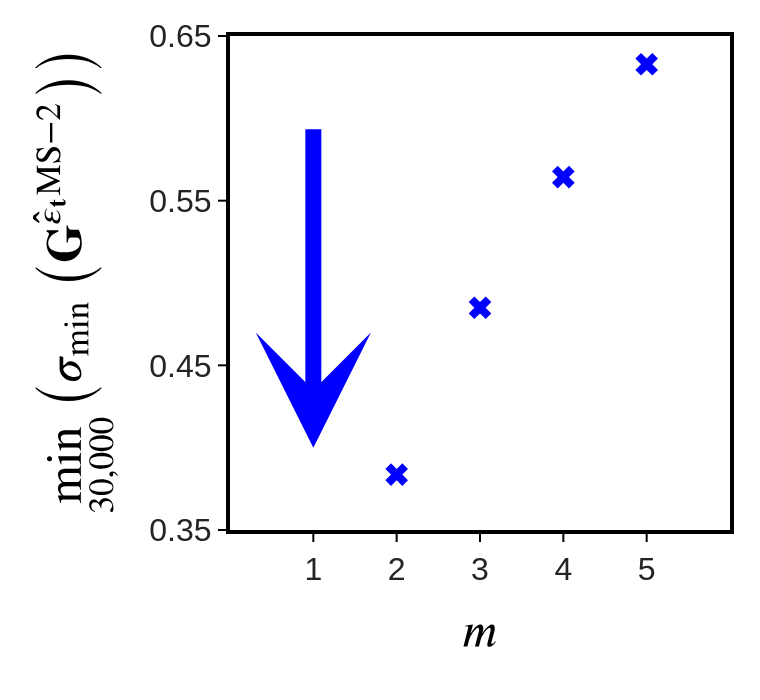}
        \caption{Minimum $\sigma_\mathrm{min}\left(\bm{G^{\hat\varepsilon_t}}\right)$ for $\delta_R=\delta_I=0.1$ for MS-2. $m=1$: one measurement at 2 GHz. $m=2,3,4,5$: $m$ measurements at equally spaced frequencies in the band 0.4-3.0 GHz, including 0.4 and 3.0 GHz. For $m=1$, $\underset{30,000}{\mathrm{min}}\left(\sigma_\mathrm{min}\left(\bm{G^{\hat\varepsilon_t}}\right)\right)$ is very low $\approx 0.01725$.} 
        \label{fig: min sigma min G MS-2}
    \end{subfigure}

    \vspace{4pt} 

    \begin{subfigure}[t]{0.42\textwidth}
        \centering
        \includegraphics[width=\linewidth]{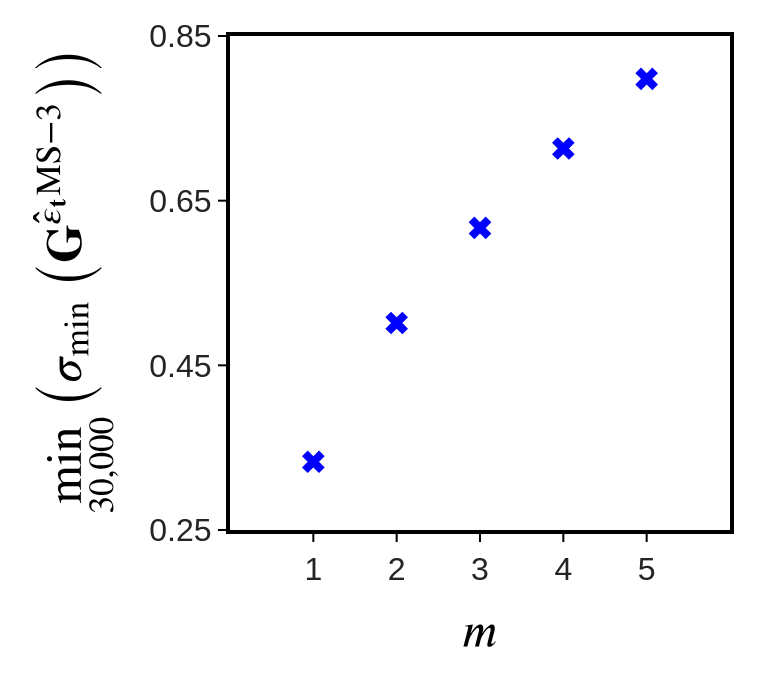}
        \caption{Minimum $\sigma_\mathrm{min}\left(\bm{G^{\hat\varepsilon_t}}\right)$ for $\delta_R=\delta_I=0.1$ for MS-3. $m=1$: one measurement at 0.5 MHz. $m=2,3,4,5$: $m$ measurements at equally spaced frequencies in the band 0.5-1.0 MHz, including 0.5 and 1.0 MHz.}
        \label{fig: min sigma min G MS-3}
    \end{subfigure}
    
    \caption{Performance of Problem \ref{eqn: Final formulation - Multifreq Epigraph} in terms of $\underset{30,000}{\mathrm{min}}\left(\sigma_\mathrm{min}\left(\bm{G^{\hat\varepsilon_t}}\right)\right)$ against number of frequencies at which measurements are taken ($m$) for $\delta_R=\delta_I=0.1$.}
    \label{fig: Problem Performance min sigma min G}
\end{figure}

The microstructure details are as follows:
\begin{itemize}
    \item Carbon-loaded epoxy material is the matrix material subscripted by `$0$'. A below percolation, moderately conductive Carbon black loaded epoxy is considered. This component is very dispersive and the Cole-Cole relaxation model with DC conductivity and with Universal Dielectric Response is used for its dispersion in the 0.5-1.0 MHz band (see \cite{cole1941dispersion}, and \cite{hill1983dielectric}):
    \begin{equation}
    \begin{gathered}
        \varepsilon_0(\omega)=\varepsilon_{0(\infty)}+\frac{\Delta\varepsilon_0}{1+(-i\omega \tau_0)^{1-\beta_0}}+\frac{i}{\omega\varepsilon_{\mathrm{vacuum}}}\left(s_{0(dc)}+A_0\omega^{s_0}\right)\text{ with,}\\
        \varepsilon_{0(\infty)} = 3.0,\> \Delta\varepsilon_0 = 150,\> \tau_0 = 2.0\times 10^{-5}\mathrm{s},\> \beta_0 = 0.35,\\
        s_{0(dc)} = 3.0\times10^{-4}\mathrm{S/m},\> A_0 = 1.0\times10^{-9}\mathrm{S/m}\cdot \mathrm{s}^{s_0},\> s_0 = 0.8
    \end{gathered}
    \end{equation}
    For $m=1$, the frequency 0.5 MHz is considered. For $m=2,3,4,5$, equally spaced frequencies in the band 0.5-1.0 MHz are considered. 
    \item Glass appears as spherical inclusions in the microstructure and is subscripted by `$1$'. The dielectric constant of glass is taken to be a constant $\varepsilon_1=5.5 + 0.05i$.
    \item Air voids or pores appear in spheres and air is subscripted by `$2$'. The dielectric constant of air is taken to be a constant $\varepsilon_2=1.0006$.
    \item $\bm{\phi^\mathrm{t}}\in\{\bm{\phi^\mathrm{t}}|\>\phi_0\ge0.65,\>\phi_1\ge0.1, \phi_2\ge0.02,\>\bm{1^3}\cdot\bm{\phi^\mathrm{t}}=1\}$
\end{itemize}
The matrix material in this material system is very dispersive and this is the best kind of material system to use with Problem \ref{eqn: Final formulation - Multifreq Epigraph} since only one frequency measurement is required and even $t^\star$ value will be low. See Figures \ref{fig: MS-1, MS-3}, \ref{fig: phi norm delta MS-3}, \ref{fig: tstar delta MS-3} and \ref{fig: min sigma min G MS-3}.

\section{Equivalent Microstructures}\label{section: Equivalent Microstructures}
\begin{figure}[H]
    \centering
    \includegraphics[width=\linewidth]{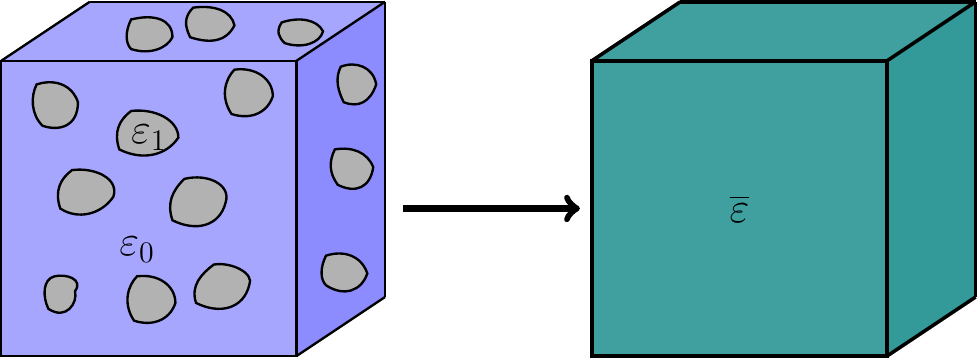}
    \caption{Forward micromechanics problem for an isotropic 2 component composite
material.}
    \label{fig: Forward Problem}
\end{figure}

\begin{figure}[H]
    \centering
    \includegraphics[width=\linewidth]{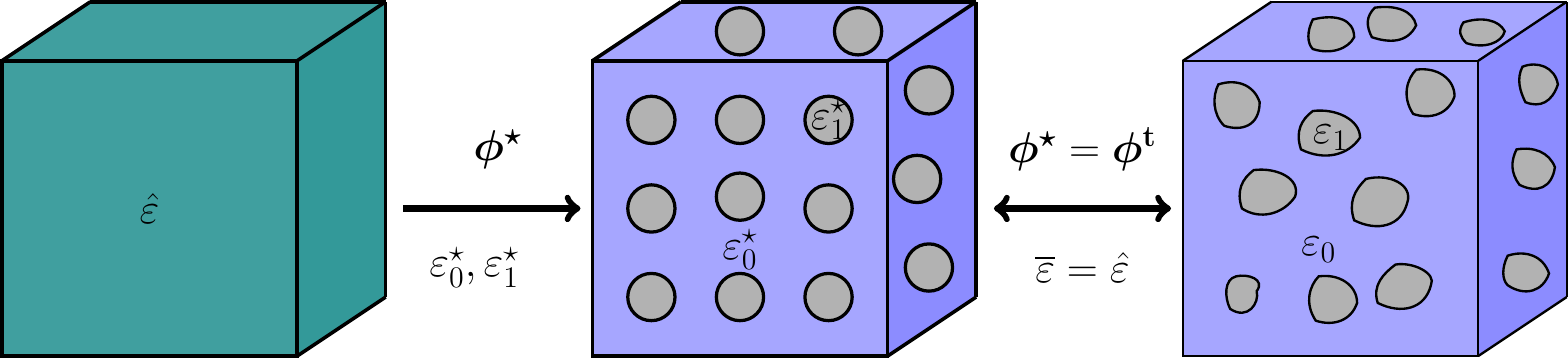}
    \caption{Inverse problem with equivalent microstructure comprising only spherical inclusions for an isotropic 2 component composite material. The goal here is to find $(\bm{\phi^\star},\varepsilon_0^\star,\varepsilon_1^\star)$ such that $\overline{\varepsilon}\left(\bm{\phi^\star},\varepsilon_0^\star,\varepsilon_1^\star\right)=\hat\varepsilon$ and $\bm{\phi^\star}=\bm{\phi^\mathrm{t}}$.
}
    \label{fig: Inverse Problem}
\end{figure}

An equivalent microstructure in this work is defined as one that has the same effective
property (dielectric constant) and volumetric compositions as the composite material under
study and has a simple microstructure of only spherical inclusions. Figures \ref{fig: Forward Problem} and \ref{fig: Inverse Problem} show the
objectives of the general forward micromechanics problem and the inverse micromechanics
problem for a 2-component isotropic composite with an equivalent microstructure comprising
only spherical inclusions. The single frequency version (no dispersion considered) of the
problem is shown here for simplicity.
The material studied must be isotropic in dielectric permittivity. It is known from micromechanics literature that such an equivalent microstructure with only spherical inclusions
with the components having the same dielectric constants as the components in the real microstructure and preserving the fractional volume split is strictly not possible for composite
materials with arbitrarily shaped inclusions (see \cite{kachanov2005quantitative}, \cite{sevostianov2012concept}, etc.). Which is why $\varepsilon_0^\star$ and $\varepsilon_1^\star$ must also be determined such that $\bm{\phi^\star}=\bm{\phi^\mathrm{t}}$. Either $\varepsilon_0^\star$ and $\varepsilon_1^\star$ can be determined some other way or an optimization problem can be formulated with $\varepsilon_0$ and $\varepsilon_1$ included as decision variables along with $\bm{\phi}$. The optimizer of the aforementioned optimization problem would be $\left(\bm{\phi^\star},\varepsilon_0^\star,\varepsilon_1^\star\right)$. For the Portland Cement Concrete material system, this must be done since the microstructure is exceedingly complex and the equivalent microstructure is very simple.

\section{Conclusion}\label{section: conclusion}
In this paper, the following objectives were met:
\begin{itemize}
    \item The inverse Eshelby-Mori-Tanaka problem for two component composites with spherical inclusions is shown to be strictly quasiconvex in $\bm\phi$.
    \item The inverse Eshelby-Mori-Tanaka problem for $n$ component composites with spherical inclusions is shown to be quasiconvex in $\bm\phi$ and the set of all minmizers is determined for $\bm\phi$ in the simplex ($\mathcal{C}$) and for $\bm\phi$ in the dominant-component ordered simplex ($\mathcal{C}_{(0)}$).
    \item An equivalent LP was formulated using the Charnes-Cooper change of variables and an epigraph formulation (\ref{eqn: epigraph with primal-dual interior pt method inclusions only}), and its inadequacies are discussed. This is done to show why either measurements at more frequencies are required or at least one complex dielectric constant measurement at a given frequency is required to solve the inverse problem for composites with 3 components.  
    \item A multi-frequency formulation of the inverse problem is determined (\ref{eqn: Final formulation - Multifreq Epigraph}). If not dielectric constant measurements at different frequencies, other information about the material system can be used to add on to the constraints of the problem. Requirements for a good prediction are $\mathrm{Rank}\left(\bm{G^{\hat\varepsilon_t}}\right)=n-1$, $\sigma_\mathrm{min}\left(\bm{G^{\hat\varepsilon_t}}\right)$ is high, and $t^\star$ is low.
    \item Problem \ref{eqn: Final formulation - Multifreq Epigraph} is then used used to validate the inverse problem formulation using forward computed values (using model \ref{eqn: MTB}) for 3 material systems. The relationship between number of frequencies at which measurements are taken, noise, and $\left\lVert\bm{\Delta\phi}\right\rVert_\infty$ is studied. Main findings at a fixed noise level are:
    \begin{itemize}
        \item The presence of a strongly dispersive component in the material system will allow excellent volumetric prediction with fewer frequency measurements.
        \item With increasing number of frequency measurements, the optimal objective $t^\star$ increases.
        \item With increasing number of frequency measurements, the lowest singular value of the constraint sensitivity matrix, $\sigma_\mathrm{min}\left(\bm{G^{\hat\varepsilon_t}}\right)$ increases.
        \item The presence of a material with strong dispersion in the imaginary part of dielectric constant (loss) will improve volumetric prediction $\left(\text{lower } \left\lVert\bm{\phi^\star}-\bm{\phi^\mathrm{t}}\right\rVert_\infty\right)$ with measurements at more frequencies ($m$).
    \end{itemize}

\end{itemize}
This paper deals with certain simple and robust optimization problems. More sophisticated problems incorporating probability density functions (encoding material symmetry other than isotropy) and other microstructural details can be formulated. The scope for work in this area of research that bridges convex optimization and Eshelby based inverse micromechanics is immense. 

\noindent\textbf{Acknowledgments:} The authors are forever indebted to Professor Kumbhako\d{n}am Ram\=ama\d{n}i R\=ajagopal for his mentoring, support, and guidance. This work would not be possible without his guidance. The first author thanks Krishna Kaushik Yanamundra and Bhaskar Vajipeyajula on helping with improvements to the paper.

\bibliographystyle{elsarticle-num} 
\bibliography{References}

@book{torquato2002random,
  title={Random heterogeneous materials: microstructure and macroscopic properties},
  author={Torquato, Salvatore and others},
  volume={16},
  year={2002},
  publisher={Springer}
}

@article{PhysRevE.59.5596,
  title = {Reconstruction of random media using Monte Carlo methods},
  author = {Manwart, C. and Hilfer, R.},
  journal = {Phys. Rev. E},
  volume = {59},
  issue = {5},
  pages = {5596--5599},
  numpages = {0},
  year = {1999},
  month = {May},
  publisher = {American Physical Society},
  doi = {10.1103/PhysRevE.59.5596},
  url = {https://link.aps.org/doi/10.1103/PhysRevE.59.5596}
}

@article{yeong1998reconstructing,
  title={Reconstructing random media},
  author={Yeong, Christofer LY and Torquato, Salvatore},
  journal={Physical review E},
  volume={57},
  number={1},
  pages={495},
  year={1998},
  publisher={APS}
}

@article{talukdar2002reconstruction,
  title={Reconstruction of chalk pore networks from 2D backscatter electron micrographs using a simulated annealing technique},
  author={Talukdar, MS and Torsaeter, O},
  journal={Journal of petroleum science and engineering},
  volume={33},
  number={4},
  pages={265--282},
  year={2002},
  publisher={Elsevier}
}

@article{talukdar2002stochastic,
  title={Stochastic reconstruction of particulate media from two-dimensional images},
  author={Talukdar, MS and Torsaeter, O and Ioannidis, MA},
  journal={Journal of colloid and interface science},
  volume={248},
  number={2},
  pages={419--428},
  year={2002},
  publisher={Elsevier}
}

@article{rintoul1997reconstruction,
  title={Reconstruction of the structure of dispersions},
  author={Rintoul, Mark D and Torquato, Salvatore},
  journal={Journal of colloid and interface science},
  volume={186},
  number={2},
  pages={467--476},
  year={1997},
  publisher={Elsevier}
}

@article{rozanski2011digital,
  title={From digital image of microstructure to the size of representative volume element: B4C/Al composite},
  author={R{\'O}{\.Z}A{\'N}SKI, ADRIAN and {\L}YD{\.Z}BA, DARIUSZ},
  journal={studia geotechnica et mechanica},
  volume={33},
  number={1},
  pages={55--68},
  year={2011}
}

@article{vcapek2009stochastic,
  title={Stochastic reconstruction of particulate media using simulated annealing: improving pore connectivity},
  author={{\v{C}}apek, P and Hejtm{\'a}nek, V and Brabec, L and Zik{\'a}nov{\'a}, A and Ko{\v{c}}i{\v{r}}{\'\i}k, M},
  journal={Transport in porous media},
  volume={76},
  pages={179--198},
  year={2009},
  publisher={Springer}
}

@article{ogierman2018inverse,
  title={Inverse identification of elastic properties of constituents of discontinuously reinforced composites},
  author={Ogierman, Witold},
  journal={Materials},
  volume={11},
  number={11},
  pages={2332},
  year={2018},
  publisher={MDPI}
}

@article{kachanov2005quantitative,
  title={On quantitative characterization of microstructures and effective properties},
  author={Kachanov, Mark and Sevostianov, Igor},
  journal={International Journal of Solids and Structures},
  volume={42},
  number={2},
  pages={309--336},
  year={2005},
  publisher={Elsevier}
}

@article{sevostianov2012concept,
  title={Is the concept of “average shape” legitimate, for a mixture of inclusions of diverse shapes?},
  author={Sevostianov, Igor and Kachanov, Mark},
  journal={International Journal of Solids and Structures},
  volume={49},
  number={23-24},
  pages={3242--3254},
  year={2012},
  publisher={Elsevier}
}

@article{Norris1989,
    author = {Norris, A. N.},
    title = "{An Examination of the Mori-Tanaka Effective Medium Approximation for Multiphase Composites}",
    journal = {Journal of Applied Mechanics},
    volume = {56},
    number = {1},
    pages = {83-88},
    year = {1989},
    month = {03},
    issn = {0021-8936},
    doi = {10.1115/1.3176070},
    url = {https://doi.org/10.1115/1.3176070},
    eprint = {https://asmedigitalcollection.asme.org/appliedmechanics/article-pdf/56/1/83/5460696/83\_1.pdf},
}

@article{BERRYMAN1996149,
title = {Critique of two explicit schemes for estimating elastic properties of multiphase composites},
journal = {Mechanics of Materials},
volume = {22},
number = {2},
pages = {149-164},
year = {1996},
issn = {0167-6636},
doi = {https://doi.org/10.1016/0167-6636(95)00035-6},
url = {https://www.sciencedirect.com/science/article/pii/0167663695000356},
author = {James G. Berryman and Patricia A. Berge}
}

@article{Esh1957,
  author = {Eshelby, J.D.},
  journal = {Proceedings of the Royal Society of London. Series A, Mathematical
	and Physical Sciences},
  number = 1226,
  pages = {376--396},
  publisher = {JSTOR},
  title = {The determination of the elastic field of an ellipsoidal inclusion,
	and related problems},
  url = {http://www.jstor.org/stable/100095},
  volume = 241,
  year = 1957
}

@book{qu2006fundamentals,
  title={Fundamentals of micromechanics of solids},
  author={Qu, Jianmin and Cherkaoui, Mohammed},
  volume={735},
address = {},
  url = {https://onlinelibrary.wiley.com/doi/pdf/10.1002/9780470117835},
  eprint = {https://onlinelibrary.wiley.com/doi/pdf/10.1002/9780470117835},
  year={2006},
  publisher={Wiley Online Library}
}

@book{mura2013micromechanics,
  title={Micromechanics of defects in solids},
  author={Mura, Toshio},
  address="Dordrecht",
  year={2013},
  publisher={Springer Science \& Business Media}
}

@article{HASHIN1963127,
title = {A variational approach to the theory of the elastic behaviour of multiphase materials},
journal = {Journal of the Mechanics and Physics of Solids},
volume = {11},
number = {2},
pages = {127-140},
year = {1963},
issn = {0022-5096},
doi = {https://doi.org/10.1016/0022-5096(63)90060-7},
url = {https://www.sciencedirect.com/science/article/pii/0022509663900607},
author = {Z. Hashin and S. Shtrikman}
}

@article{lydzba2019principle,
  title={Principle of equivalent microstructure in micromechanics and its connection with the replacement relations. Thermal conductivity problem},
  author={{\L}yd{\.z}ba, D and R{\'o}{\.z}a{\'n}ski, A and Sevostianov, I and Stefaniuk, D},
  journal={International Journal of Engineering Science},
  volume={144},
  pages={103126},
  year={2019},
  publisher={Elsevier}
}

@article{LYDZBA201820,
title = {Equivalent microstructure problem: Mathematical formulation and numerical solution},
journal = {International Journal of Engineering Science},
volume = {123},
pages = {20-35},
year = {2018},
issn = {0020-7225},
doi = {https://doi.org/10.1016/j.ijengsci.2017.11.007},
url = {https://www.sciencedirect.com/science/article/pii/S0020722517318980},
author = {Dariusz Łydżba and Adrian Różański and Damian Stefaniuk}
}

@INPROCEEDINGS{ECOS,

  author={Domahidi, Alexander and Chu, Eric and Boyd, Stephen},

  booktitle={2013 European Control Conference (ECC)}, 

  title={ECOS: An SOCP solver for embedded systems}, 

  year={2013},

  volume={},

  number={},

  pages={3071-3076},

  doi={10.23919/ECC.2013.6669541}}

@book{boyd2004convex,
  title={Convex optimization},
  author={Boyd, Stephen and Vandenberghe, Lieven},
  year={2004},
  publisher={Cambridge university press}
}

@article{benveniste1986effective,
  title={On the effective thermal conductivity of multiphase composites},
  author={Benveniste, Ya},
  journal={Zeitschrift f{\"u}r angewandte Mathematik und Physik ZAMP},
  volume={37},
  number={5},
  pages={696--713},
  year={1986},
  publisher={Springer}
}

@book{maxwell1873treatise,
  title={A treatise on electricity and magnetism},
  author={Maxwell, James Clerk},
  volume={1},
  year={1873},
  publisher={Clarendon press}
}

@book{stratton2007electromagnetic,
  title={Electromagnetic theory},
  author={Stratton, Julius Adams},
  year={2007},
  publisher={John Wiley \& Sons}
}

@article{debye1929polar,
  title={Polar molecules, the chemical catalog company},
  author={Debye, P},
  journal={Inc., New York},
  volume={89},
  year={1929}
}

@article{cole1941dispersion,
  title={Dispersion and absorption in dielectrics I. Alternating current characteristics},
  author={Cole, Kenneth S and Cole, Robert H},
  journal={The Journal of chemical physics},
  volume={9},
  number={4},
  pages={341--351},
  year={1941},
  publisher={American Institute of Physics}
}

@article{hill1983dielectric,
  title={The dielectric behaviour of condensed matter and its many-body interpretation},
  author={Hill, RM and Jonscher, AK},
  journal={Contemporary Physics},
  volume={24},
  number={1},
  pages={75--110},
  year={1983},
  publisher={Taylor \& Francis}
}

@article{thaler2013bounds,
  title={Bounds on the volume of an inclusion in a body from a complex conductivity measurement},
  author={Thaler, Andrew E and Milton, Graeme W},
  journal={arXiv preprint arXiv:1306.6608},
  year={2013}
}

@book{cherkaev2012variational,
  title={Variational methods for structural optimization},
  author={Cherkaev, Andrej},
  volume={140},
  year={2012},
  publisher={Springer Science \& Business Media}
}

@article{milton2012bounds,
  title={Bounds on the volume fraction of 2-phase, 2-dimensional elastic bodies and on (stress, strain) pairs in composites},
  author={Milton, Graeme Walter and Nguyen, Loc Hoang},
  journal={Comptes Rendus M{\'e}canique},
  volume={340},
  number={4-5},
  pages={193--204},
  year={2012},
  publisher={Elsevier}
}

@article{kang2013bounds,
  title={Bounds on the volume fractions of two materials in a three-dimensional body from boundary measurements by the translation method},
  author={Kang, Hyeonbae and Milton, Graeme W},
  journal={SIAM Journal on Applied Mathematics},
  volume={73},
  number={1},
  pages={475--492},
  year={2013},
  publisher={SIAM}
}

@article{hashin1983analysis,
    author = {Hashin, Z.},
    title = "{Analysis of Composite Materials—A Survey}",
    journal = {Journal of Applied Mechanics},
    volume = {50},
    number = {3},
    pages = {481-505},
    year = {1983},
    month = {09},
    issn = {0021-8936},
    doi = {10.1115/1.3167081},
    url = {https://doi.org/10.1115/1.3167081}
}

\end{document}